\documentclass[a4paper,dvipsnames,fleqn]{article}

\usepackage{amsmath}
\usepackage{amssymb}
\usepackage{amsthm}
\allowdisplaybreaks
\usepackage{bm}
\usepackage{enumerate}
\usepackage[backend=biber,bibstyle=alphabetic,citestyle=alphabetic,sorting=nyt,maxnames=6,giveninits=false]{biblatex}
	\AtEveryBibitem{\clearfield{doi}}
	\AtEveryBibitem{\clearfield{isbn}}
	\AtEveryBibitem{\clearfield{issn}}
	\AtEveryBibitem{\clearfield{url}}
	\renewbibmacro{in:}{\ifentrytype{article}{}{\printtext{\bibstring{in}\intitlepunct}}}
	\AtEveryBibitem{\clearfield{number}}
	\AtEveryBibitem{\ifentrytype{book}{\clearfield{pages}}{}}
	\DeclareFieldFormat[article,incollection,inbook,unpublished]{title}{#1}
	\DeclareFieldFormat[article]{volume}{\textrm{#1}}
	\DeclareFieldFormat[article]{pages}{#1}
	\DeclareFieldFormat[book]{edition}{#1}
\usepackage[colorlinks=true,allcolors=purple]{hyperref}
\usepackage{cleveref}
	\crefformat{equation}{#2(#1)#3}
\usepackage[british]{babel}
\usepackage{csquotes}
\usepackage[UKenglish]{isodate}
\usepackage{dsfont}
\usepackage{enumitem}
\usepackage[OT2,OT1]{fontenc}
\usepackage[cal=cm,scr=boondoxo]{mathalfa}
\usepackage{pifont}
\usepackage{stmaryrd}
\usepackage{textcomp}
\usepackage{tikz}
	\usetikzlibrary{arrows}
	\usetikzlibrary{patterns}
\usepackage{tikz-cd}
\usepackage[textwidth=3cm,colorinlistoftodos]{todonotes}
\usepackage{xfrac}
\usepackage{csquotes}
\AtBeginEnvironment{quote}{\itshape}
\numberwithin{equation}{section}			% Enumeration of equations
\swapnumbers						% Enumeration of theorem-like environments

\newcommand\cyr{%					% cyrillic font
\renewcommand\rmdefault{wncyr}%
\renewcommand\sfdefault{wncyss}%
\renewcommand\encodingdefault{OT2}%
\normalfont
\selectfont}
\DeclareTextFontCommand{\textcyr}{\cyr}

\newcounter{Enum}				% Enumerated list
\newenvironment{Enumerate}{\begin{enumerate}[label={\rm({\roman*})}]}{\end{enumerate}}

\newcommand{\descriptionlabelsave}{}		% Itemized list
\newenvironment{Itemize}{%
	\renewcommand{\descriptionlabelsave}{\descriptionlabel}\renewcommand{\descriptionlabel}{$\triangleright$}%
	\begin{description}[leftmargin=15pt,itemindent=-5.2pt]}{%
	\end{description}\renewcommand{\descriptionlabel}{\descriptionlabelsave}}

\newcounter{StepsCount}				% Enumerated list with no indentation (e.g. steps in proof)

\newcounter{StepsRefCount}

\theoremstyle{plain}
	\newtheorem{lemma}{Lemma}[section]
	\newtheorem{proposition}[lemma]{Proposition}
	\newtheorem{theorem}[lemma]{Theorem}
	\newtheorem{corollary}[lemma]{Corollary}
	\newcommand{\GenericTheoremName}{}\newtheorem{generictheorem}[lemma]{\GenericTheoremName}
\theoremstyle{definition}
	\newtheorem{definition}[lemma]{Definition}
	\newcommand{\GenericDefinitionName}{}\newtheorem{genericdefinition}[lemma]{\GenericDefinitionName}
	\newtheorem{notation}[lemma]{Notation}
\theoremstyle{remark}
	\newtheorem{remark}[lemma]{Remark}
	\newtheorem{example}[lemma]{Example}
	\newcommand{\GenericRemarkName}{}\newtheorem{genericremark}[lemma]{\GenericRemarkName}
\newenvironment{GenericTheorem}[1]
	{\renewcommand{\GenericTheoremName}{#1}\begin{generictheorem}}{\par\noindent\centerline{\rule{5em}{1pt}}\end{generictheorem}}%6.5
\newcommand{\mc}[1]{{\mathcal{#1}}}			% --- abbreviation ---
\newcommand{\ms}[1]{{\mathscr{#1}}}			% --- abbreviation ---
\newcommand{\bb}[1]{{\mathbb{#1}}}			% --- abbreviation ---
\newcommand{\wt}{\widetilde}				% --- abbreviation ---
\DeclareMathOperator{\tr}{tr}

\newcommand{\Side}[1]{\hfill{#1}\kern10pt}		% text put on the right side of line with offset
\newcommand{\FD}[5]{%					% definition of function from {#1} to {#2} mapping {#4} to {#5}
	\DF\left\{\begin{array}{rcl}{#1}&\to &{#2}\\[#3pt] {#4}&\mapsto &{#5}\end{array}\right.}

\newcommand{\smmatrix}[4]{\Bigl(			% small matrix for use in textline
\begin{smallmatrix}
\hspace*{-0.2ex} #1 \hspace*{0.2ex} & \hspace*{0.2ex} #2 \hspace*{-0.2ex}
\\[0.5ex]
\hspace*{-0.2ex} #3 \hspace*{0.2ex} & \hspace*{0.2ex} #4 \hspace*{-0.2ex}
\end{smallmatrix}
\Bigr)}
\newcommand{\Dummy}{\text{\textvisiblespace\kern1pt}}	% Platzhaltersymbol fuer Funktionsargumente
\newcommand{\Smallo}{{\rm o}}				% small o
\newcommand{\BigO}{{\rm O}}				% big o

\newcommand{\DQ}{\mkern6mu}				% distance for successive quantors
\newcommand{\DP}{{.\kern5pt}}				% delimiter for predicate formula
\newcommand{\DF}{\colon}				% delimiter for function domain/codomain
\newcommand{\DE}{\mathrel{\mathop:}=}			% defining equality
\newcommand{\DD}{\mkern4mu\mathrm{d}}			% distance and rm-d for integration differential

\begin{document}

\begin{flushleft}
	{\Large\bf Nevanlinna matrix estimates without regularity conditions}
	\\[5mm]
	\textsc{
	Jakob Reiffenstein
	\hspace*{-14pt}
		\renewcommand{\thefootnote}{\fnsymbol{footnote}}
		\setcounter{footnote}{2}
			\footnote{
The author was supported by the project I 4600 of the Austrian Science Fund
(FWF), and by the Sverker Lerheden
foundation.}
	} \\[1ex]
	{\small
	\textbf{Abstract.}
The Nevanlinna matrix of a half-line Jacobi operator coincides, up to multiplication with a constant matrix, with the monodromy matrix of an associated canonical system. This canonical system is discrete in a certain sense, and is determined by two sequences, called ``lengths'' and ``angles''. We derive new lower and upper estimates for the norm of the monodromy matrix in terms of the lengths and angles, without imposing any restrictions on these sequences. Returning to the Jacobi setting, we show that the order of the Nevanlinna matrix is always greater than or equal to the convergence exponent of the off-diagonal sequence of Jacobi parameters, and obtain a generalisation of a classical theorem of Berezanskii.
	\\[3mm]
	}
\end{flushleft}

%%%%%%%%%%%%%%%%%%%%%%%%%%%%%%%%%%%%%%%%%%%%%%%%%%%%%%%%%%%%%%%%%%%%%%%%%%%%%%%%%%%%%%%%%
% disable in final version
%\pagenumbering{roman}
%\fbox{
%\parbox{100mm}{
%\hspace*{0pt}\\[1mm]
%\centerline{{\Large\ding{45}}\quad\,{\large\sc Draft}%\quad{\Large\ding{45}}}
%\hspace*{0pt}\\[-2mm]
%\textcircledP\ \ Preliminary version Fri 8 Jul 2023 16:46
%\\[2mm]
%\cleanlookdateon
%\hspace*{5mm} Compilation date: \today
%\\[2mm]
%\ding{233}\quad Use pdflatex/biber to compile
%\\[-1mm]
%}
%}
%\tableofcontents
%\listoftodos
%\newpage
%\pagenumbering{arabic}
%\setcounter{page}{1}
%
%%%%%%%%%%%%%%%%%%%%%%%%%%%%%%%%%%%%%%%%%%%%%%%%%%%%%%%%%%%%%%%%%%%%%%%%%%%%%%%%%%%%%%%%%

%---------
%   TEXTBODY
%---------

%\newpage

%
%
%
%\section[\textcolor{ForestGreen}{}]{}
%\section[\textcolor{Dandelion}{}]{}
%\section[\textcolor{BrickRed}{...}]{...}

\section{Introduction}

We investigate two equivalent types of operators. On the one hand, we consider two-dimensional canonical systems
\begin{align}
\label{B36}
y'(t)=zJH(t)y(t), \qquad t \in (0,L) \text{ a.e.},
\end{align}
where $L \in (0,\infty ]$, $J:=\smmatrix 0{-1}10$, $z \in \mathbb{C}$, and the \emph{Hamiltonian} $H$ is a measurable $2 \times 2$ matrix-valued function on $(0,L)$ that satisfies 
\[
H(t) \geq 0, \quad \tr H(t)=1, \qquad t \in (0,L) \text{ a.e.}
\]
We assume that $H$ is a \emph{Hamburger Hamiltonian}, i.e., it is of a particular discrete form defined further below, cf. \Cref{S21}.

On the other hand, we investigate unbounded Jacobi matrices given by
\begin{align}
\label{K71}
\ms J=\begin{pmatrix}
a_0 & b_0 & 0 & & & \\
b_0 & a_1 & b_1 & 0 & &\\
0   & b_1 & a_2 & b_2 & & \\
    & \raisebox{5pt}[0pt][5pt]{$0$}   & \raisebox{5pt}[0pt][5pt]{$b_2$} & \ddots & \ddots & \\
    &     &     & \ddots & \ddots & \\
\end{pmatrix}
\end{align}
with $a_n \in \mathbb{R}$ and $b_n>0$. There is a natural way of rewriting the formal eigenvalue equation $\ms Ju=zu$ into \eqref{B36}, with the Hamiltonian $H$ being a Hamburger Hamiltonian. This transformation is injective and such that the minimal operator of $\ms J$ in $\ell^2(\bb N)$ is unitarily equivalent to a certain one-dimensional extension of the minimal operator associated with $H$. For this reason we think of $\ms J$ and $H$ interchangeably.

The minimal operator of $\ms J$ is either self-adjoint (limit point case) or symmetric with deficiency indices $(1,1)$ (limit circle case). In Hamiltonian parameters, the limit circle case occurs if and only if $L<\infty$, while in terms of the Jacobi parameters $a_n$ and $b_n$ there is no simple condition. We always assume that limit circle case takes place, i.e., $L<\infty$. 

Since we only consider limit circle case, the spectrum of any self-adjoint extension of the minimal operator of $\ms J$ is discrete. There exists a $2 \times 2$ matrix-valued entire function $W$, called \emph{Nevanlinna matrix}, such that the spectrum of any self-adjoint extension coincides with the set of zeros of a linear combination of the entries of $W$. An analogous object exists for canonical systems, i.e., for every Hamburger Hamiltonian $H$ there is a matrix function $W_H$, called \emph{monodromy matrix}, with the properties described above. When $\ms J$ and $H$ are related by unitary equivalence of the operators mentioned above, $W$ and $W_H$ coincide up to constant multiplicative factors. Thus we refer to $W$ and $W_H$ interchangeably, too.  \\

%In this paper we estimate the growth of the monodromy matrix for canonical systems \eqref{B36} when $H$ is a Hamburger Hamiltonian, also yielding new estimates for the Nevanlinna matrix of a Jacobi operator directly in terms of $a_n,b_n$. Thanks to the connection between the growth of an entire function and the distribution of its zeros, we obtain estimates for the density of eigenvalues of self-adjoint realisations of \eqref{B36}. \\
All four entries of $W$ have similar growth, in the sense that the quotient of any two entries is (up to a possible minus sign) a product of at most two Herglotz functions \cite[Theorem 4.18]{remling:2018}. Since Herglotz functions cannot grow or decay too fast, one can show that the entries $w_{ij}$ of $W$ all have the same order
\begin{align*}
\rho &\DE \limsup_{r \to \infty} \frac{\log \log \max_{|z|=r} |w_{ij}(z)|}{\log r}, \quad i,j=1,2.
\end{align*}
According to a theorem of M. Riesz \cite{riesz:1923a,akhiezer:1961} (or, alternatively, the Krein-de~Branges formula \cite{krein:1951,debranges:1961}), the matrix function $W$ is of minimal exponential type. In particular this implies that $\rho \in [0,1]$. 

Let us give an instance of how the growth of $W$ is related to the density of the eigenvalues of a self-adjoint realisation $A$ of the Jacobi matrix or canonical system.
Choosing an enumeration $(\lambda_n)_{n \in \bb N}$ of the nonzero eigenvalues of $A$ we are interested in the convergence exponent
\begin{align}
\label{K76}
\mc E ((\lambda_n)_{n \in \bb N}) \DE\inf \big\{\alpha>0 \DF  \sum_{n \in \bb N} |\lambda_n|^{-\alpha} <\infty \big\},
\end{align}
which is a measure for how dense the eigenvalues are. This number is independent of $A$, since any two self-adjoint realisations are perturbations of rank at most two of each other. 
%Hence, we can assume without loss of generality that $A$ is the self-adjoint realisation whose spectrum coincides with the set of zeros of $w_{22}$. 
It follows from the classical theory of entire functions that
\begin{align}
\label{K73}
\mc E ((\lambda_n)_{n \in \bb N})=\rho,
\end{align}
in other words, knowledge about the growth of $W$ translates to knowledge about the density of eigenvalues of any self-adjoint realisation of \eqref{B36} (and the equality \eqref{K73} is just one instance of this connection).

The problem of estimating the growth of the monodromy or Nevanlinna matrix (e.g., by estimating $\rho$) has been treated in a fair number of papers \cite{livshits:1939,berezanskii:1956,berg.szwarc:2014,romanov:2017,pruckner.romanov.woracek:jaco,pruckner.woracek:srt,pruckner:blubb,pruckner.reiffenstein.woracek:sinqB-arXiv,reiffenstein:kac-hamA-arXiv}. While these papers feature a large variety of upper estimates, the go-to lower bound for $\rho$ in almost all of them is $\mc E \big((b_n)_{n=0}^\infty \big)$, i.e., the convergence exponent of the off-diagonal sequence $(b_n)_{n=0}^\infty$. It is equal to $\rho$ in surprisingly many cases. For instance, Berezanskii's Theorem \cite{berezanskii:1956,berg.szwarc:2014} gives sufficient conditions, in terms of $a_n,b_n$, for $\ms J$ to be in limit circle case and for $\rho =\mc E ((b_n)_{n=0}^\infty )$ to hold.
%In fact, there might be only one known example of a limit circle Jacobi matrix for which the order of $W$ can be computed and is different from the convergence exponent of $(b_n)_{n=0}^\infty$, cf. \cite[Theorem~1.3]{reiffenstein:kac-hamA-arXiv}. 
However, previous proofs of the lower estimate $\rho \geq \mc E ((b_n)_{n=0}^\infty )$ all required regularity conditions on the data, e.g., log-concavity or regular variation of $(b_n)_{n=0}^\infty$.

In this paper we derive lower and upper bounds for $\rho$ in terms of the parameter sequences of $\ms J$ and $H$, whose growth or decay we measure in terms of convergence exponents. Our results largely do not depend on regularity conditions, since convergence exponents are always well-defined. We present all lower bounds in \Cref{S3} and all upper bounds in \Cref{S53}.

Regarding lower bounds we would like to highlight \Cref{K65} which, in particular, states that $\rho \geq \mc E ((b_n)_{n=0}^\infty )$ for any Jacobi matrix in limit circle case. This result is the most illustrative but roughest out of a number of lower estimates for the Nevanlinna matrix featured in \Cref{S52}; for instance, we also give a pointwise lower estimate for $W$ along the imaginary axis.

\Cref{S53} is dedicated to upper bounds, mostly in terms of the Hamiltonian $H$. The foundational theorem of that section is \Cref{K26}, where we give an upper bound for an entry of $W_H$ assuming we can estimate the determinants of the integrals of $H$ over subintervals of $[0,L]$ in a specific way. Different methods of estimating these determinants result in different upper bounds for $W_H$. Which of these methods yields the best results depends mostly on one sequence of parameters of $H$, which we call \textit{angles} (we give a definition in \Cref{S21}). In Sections \ref{S533} and \ref{S532} we prove the following upper bounds for (an entry of) $W_H$, each of which is best suited for a specific behaviour of the angles:
\begin{itemize}
\item[$\rhd$] \Cref{K89} and \Cref{K66} are applicable when the differences of consecutive angles are summable, and they give upper bounds that often coincide with the universal lower bound from \Cref{K65} on the level of order.
\item[$\rhd$] \Cref{K79} is applicable to any Hamburger Hamiltonian and best suited for convergent angles, measuring the speed of convergence in $\ell^p$-sense. It gives a good upper bound if the angles converge quickly.
\item[$\rhd$] \Cref{K49} is also applicable to any Hamburger Hamiltonian. It factors in the decay of the differences of consecutive angles (in the $\ell^p$-sense) as well as convergence of the angles (but not anymore in $\ell^p$-sense). 
\end{itemize}
We illustrate these results in \Cref{S534} by means of an example featuring mixed rates of decay of the sequence of \emph{lengths} of $H$ (see again \Cref{S21} for a definition). 

We conclude the paper by returning to the setting of Jacobi matrices. In \Cref{K84} we give a simple condition on the growth of the solutions of the Jacobi recurrence and show that it implies that limit circle case takes place and that the order of $W$ is equal to the convergence exponent of $(b_n)_{n=0}^\infty$. We put this condition into context in \Cref{K85}, generalising Berezanskii's theorem by considerably relaxing the classical assumptions while keeping everything explicit in terms of $a_n,b_n$.

\section{Preliminaries}

We introduce Hamburger Hamiltonians, which are the central objects considered in this paper, and give a short collection of notions and results that we are going to use extensively.

\subsection{Hamburger Hamiltonians}
\label{S21}
Let $(l_j)_{j=1}^\infty$ be a sequence of positive numbers and $(\phi_j)_{j=1}^\infty$ be a sequence of real
	numbers satisfying $\phi_{j+1}-\phi_j \not\equiv 0 \mod \pi$. Using the notation
	\begin{align*}
	x_0 &\DE 0, \qquad x_n \DE\sum_{j=1}^n l_j, \qquad n \in \mathbb{N} , \\  L &\DE\sum_{j=1}^\infty l_j
	\end{align*}
	and
\[
\xi_\phi \DE \binom{\cos \phi}{\sin \phi},
\]
	we define a Hamiltonian $H$ on $(0,L)$ by setting 
	\begin{align}
		H(t)\DE\xi_{\phi_j}\xi_{\phi_j}^\top\quad\text{for }j\in\bb N\text{ and }
		x_{j-1} \leq t< x_j
		.
	\end{align}
	\begin{center}
	\begin{tikzpicture}[x=1.2pt,y=1.2pt,scale=0.8,font=\fontsize{8}{8}]
		\draw[thick] (10,30)--(215,30);
		\draw[dotted, thick] (215,30)--(270,30);
		\draw[thick] (10,25)--(10,35);
		\draw[thick] (70,25)--(70,35);
		\draw[thick] (120,25)--(120,35);
		\draw[thick] (160,25)--(160,35);
		\draw[thick] (190,25)--(190,35);
		\draw[thick] (210,25)--(210,35);
		\draw[thick] (270,25)--(270,35);
		\draw (40,44) node {${\displaystyle \xi_{\phi_1}\xi_{\phi_1}^\top}$};
		\draw (95,44) node {${\displaystyle \xi_{\phi_2}\xi_{\phi_2}^\top}$};
		\draw (140,44) node {${\displaystyle \xi_{\phi_3}\xi_{\phi_3}^\top}$};
		\draw (177,43) node {${\cdots}$};
		\draw (-20,30) node {\large $H\!:$};
		\draw (10,18) node {${\displaystyle x_0}$};
		\draw[dashed,stealth-stealth] (11,26)--(69,26);
		\draw (40,21) node {${l_1}$};
		\draw (70,18) node {${x_1}$};
		\draw[dashed,stealth-stealth] (71,26)--(119,26);
		\draw (95,21) node {${l_2}$};
		\draw (120,18) node {${x_2}$};
		\draw[dashed,stealth-stealth] (121,26)--(159,26);
		\draw (140,21) node {${l_3}$};
		\draw (160,18) node {${x_3}$};
		\draw (195,18) node {${\cdots}$};
		\draw (270,18) node {${\displaystyle L}$};
	\end{tikzpicture}
	\end{center}
	
	The numbers $l_j$ and $\phi_j$ are referred to as \emph{lengths} and \emph{angles} of the Hamiltonian.

We only consider Hamburger Hamiltonians satisfying $L<\infty$. As was mentioned earlier, this is equivalent to the prevalence of limit circle case. The condition $\phi_{j+1}-\phi_j \not\equiv 0 \mod \pi$ makes sure that $H$ has a jump at every $x_j$, and hence the lengths are unique.
\newline

%Hence, when we make an assumption on ``the sequence of angles'' of a Hamburger Hamiltonian $H$, what we mean is that \emph{there exists} $(\phi_j)_{j=1}^\infty$ such that $H=H_{l,\phi}$ and that satisfies the .

%Given a Hamburger Hamiltonian $H_{l,\phi}$, we will sometimes ask for $(\phi_j)_{j=1}^\infty$ to satisfy some condition $(\ast )$. We point out that even if $(\ast )$ does not hold for $(\phi_j)_{j=1}^\infty$, it might be satisfied by another sequence $(\tilde \phi_j)_{j=1}^\infty$ for which $H_{l,\phi}=H_{l,\tilde \phi}$.

%This means that 
%In order to avoid confusion, we point out that for a given Hamburger Hamiltonian $H_{l,\phi}$ and a given condition $(\ast )$ to be satisfied
	
\noindent Let $W_H=(w_{H,ij})_{i,j=1}^2$ be the fundamental solution of (the transpose of) \eqref{B36}:
\begin{align*}
\begin{cases}
\frac{d}{dt} W_H(t;z)J=zW_H(t;z)H(t), &t \in (0,L), \\
W_H(0,z)=I.
\end{cases}
\end{align*}
Since limit circle case takes place at $L$, the solution $W_H$ exists up to $t=L$, and we define the monodromy matrix as $W_H(z)\DE W_H(L;z)$.

\subsection{Comparing two functions to each other}
Consider functions $f,g: X \to (0,\infty)$, where $X$ is any set.
We use the following notation to compare $f$ and $g$.

\begin{itemize}
\item[$\rhd$] We write $f \lesssim g$ (or $f(x) \lesssim g(x)$) if there exists $C>0$ such that $f(x) \leq Cg(x)$ for all $x \in X$. By $f \gtrsim g$ we mean $g \lesssim f$, and $f \asymp g$ is used for ``$f \lesssim g$ and $f \gtrsim g$''.
\item[$\rhd$] Let $Y \subseteq X$. We write ``$f \lesssim g$ on $Y$'' if $f|_Y \lesssim g|_Y$.
\item[$\rhd$] If $X=[1,\infty)$ or $X=\bb N$ (or any other directed set) we write ``$f(x) \lesssim g(x)$ for sufficiently large $x$'' to say that there exists $x_0 \in X$ such that $f \lesssim g$ 
on $\{x \in X \DF \, x \geq x_0 \}$. 
\item[$\rhd$] The notation $f=\BigO (g)$ stands for $\limsup_{x \to \infty} \frac{f(x)}{g(x)}<\infty$.
\end{itemize}

\subsection{A method of estimating \boldmath{$W_H$}}

%\begin{definition}
%Let $H$ be a Hamiltonian on $(0,L)$, and set $r_0 \DE \frac{1}{\sqrt{\det \Omega_H(0,L)}}$. For $r \geq r_0$, let $\hat t(r) \in (0,L]$ be defined implicitly by the equation
%\begin{align*}
%\det \Omega_H(0,\hat t(r))=\frac{1}{r^2},
%\end{align*}
%which has a unique solution due to the properties of $\det \Omega_H$ stated in \Cref{K40}. \\
%If $r \geq r_0$ and $t \in [\hat t(r),L]$, let $\hat s(t;r)$ be defined implicitly by
%\begin{align*}
%\det \Omega_H(\hat s(t;r),t)=\frac{1}{r^2},
%\end{align*}
%which is again possible due to \Cref{K40}.
%\end{definition}

When measuring the growth of $W_H$ we usually examine its lower right entry $w_{H,22}$ and its growth along the imaginary axis. In fact, $|w_{H,22}(ir)|=\max_{|z|=r} |w_{H,22}(z)|$ and hence we can calculate $\rho$ given $|w_{H,22}(ir)|$ only. Here the claimed equality follows from the product representation 
\begin{align*}
w_{H,22}(z)=\lim_{R \to \infty} \prod_{|\lambda_n|<R} \Big(1-\frac{z}{\lambda_n} \Big),
\end{align*}
which, up to a multiplicative constant, is true for any entire function of exponential type that is real along the real axis and of bounded type in the upper half-plane \cite[V.Theorem~11]{levin:1980}.

In \cite[Definition~3.1]{langer.reiffenstein.woracek:kacest-arXiv} a function $K_H$ is defined, which is derived explicitly from $H$ and satisfies
	\begin{equation}\label{X65}
		\log |w_{H,22}(ir)| \asymp \int_0^{L} K_H(t;r)\DD t,
		\qquad r>0.
	\end{equation}
The precise form of $K_H$ is not relevant for our purposes, as we will only need more explicit estimates in the form of \Cref{K2,X67} below. Both the lower and upper constant implicit in $\asymp$ in \eqref{X65} are universal, i.e., independent of $H$. 

Another important role is played by the function
\begin{align}
\Omega (s,t)=\int_s^t H(x) \, dx.
\end{align}
We can express $\det \Omega$ more explicitly in terms of the Hamiltonian parameters:
\begin{align}
\label{KX54}
\det\Omega(x_m,x_n) =\frac 12 \sum_{j,k=m+1}^n l_jl_k \sin^2(\phi_j-\phi_k).
\end{align}
This formula is a special case of \cite[Lemma~6.3]{langer.reiffenstein.woracek:kacest-arXiv}.

We now recall two central tools for estimating $W_H$ from above and below.

\begin{lemma}[{\cite[Lemma~5.5]{langer.reiffenstein.woracek:kacest-arXiv}}]
\label{K2}
Let $r>0$, and assume we have points $s_0,s_1,s_2$ such that
\begin{align*}
0 \leq s_0 < s_1 < s_2 \leq L, \qquad \det \Omega (s_0,s_1) \geq \frac{1}{r^2} \geq \det \Omega (s_1,s_2).
\end{align*}
Then
\begin{align}
\label{K24}
	 \int_{s_1}^{s_2}  K_H(t;r) \DD t \le e \log \big(r^2\det\Omega(s_0,s_2)\big).
\end{align}
\end{lemma}

\begin{lemma}[{\cite[Lemma~5.7]{langer.reiffenstein.woracek:kacest-arXiv}}]
\label{X67}
Let $r>0$. Assume we have points $0 \leq s_0<s_1<\ldots <s_k \leq L$ such that
	\begin{equation}\label{X61}
		\det\Omega(s_{j-1},s_j) \ge \frac{1}{r^2}, \qquad j\in\{1,\ldots,k\}.
	\end{equation}
	Then we have
	\[
		\log |w_{H,22}(ir)| \geq c \big(
		k\log 2-\log (Lr) \big),
	\]
	where $c$ is a universal constant, i.e., it does not depend on the Hamiltonian or on the points $s_0,\ldots,s_k$.
\end{lemma}

If we only care about large $r$ and are willing to accept dependence of the constant on the Hamiltonian, we can neglect the additive term $\log (Lr)$. In fact, since $w_{H,22}(x_2;\cdot)$ is a second-degree polynomial, we have 
\[
\log |w_{H,22}(ir)| \geq \log |w_{H,22}(x_2;ir)| \geq c'\log (Lr)>0
\]
for sufficiently large $r$. Consequently,
\begin{align}
\label{K101}
\begin{split}
\log |w_{H,22}(ir)| &\geq \min\{c',c\} \max \{\log (Lr),k\log 2-\log (Lr) \} \\
&\geq \min\{c',c\} \frac{k\log 2}{2}=c'' \cdot k,
\end{split}
\end{align}
where $c''$ depends on $H$ but not on $k$.

\subsection{The relation between Hamiltonian and Jacobi parameters}

Let $H$ be a Hamburger Hamiltonian with lengths $(l_j)_{j=1}^\infty$ and angles $(\phi_j)_{j=1}^\infty$. Some of the results that we obtain in terms of the lengths and angles can be transferred to Jacobi parameters immediately, thanks to the explicit relation between Hamiltonian and Jacobi parameters. 
Namely, if a Jacobi matrix $\ms J$ is given, we define the corresponding Hamburger Hamiltonian $H$ by setting $l_1 \DE 1$, $\phi_1 \DE \pi /2$, and finding the remaining lengths and angles by recursively solving
\begin{align}
\label{K11}
a_0 &= \tan \phi_2, \\
\label{K44}
a_n &= - \frac{\sin (\phi_{n+2} - \phi_{n})}{l_{n+1} \sin (\phi_{n+2} -\phi_{n+1} ) \sin(\phi_{n+1} - \phi_{n})}, &&n=1,2,\ldots , \\
\label{K62}
b_n &= \frac{1}{\sqrt{ l_{n+1}l_{n+2}} |\sin (\phi_{n+2} - \phi_{n+1})|}, &&n =  0,1 ,2 ,\ldots.
\end{align}
Note that the angles $\phi_n$ are only determined modulo $\pi$, but this still defines a Hamburger Hamiltonian uniquely. When the parameters of $\ms J$ and $H$ are related in this way, the minimal operator of $\ms J$ is unitarily equivalent to the operator associated with \eqref{B36} and boundary conditions $y_1(0)=0$, $y_1(L)=y_2(L)=0$. For more detailed information see \cite{kac:1999}.

The following relation is then true for the Nevanlinna matrix $W$ of $\ms J$ and the monodromy matrix $W_H$ of $H$:
\[
\begin{pmatrix}
A(z) & -C(z) \\
-B(z) & D(z)
\end{pmatrix}
=
\begin{pmatrix}
1 & 0 \\
0 & -1
\end{pmatrix}
W(z)
\begin{pmatrix}
1 & 0 \\
0 & -1
\end{pmatrix}
=
W_H(z)J.
\]
In particular, we have $w_{H,22}=-B$.

\section{Universal lower bounds}
\label{S3}

The aim of this section is to present lower bounds for the growth of $W_H$ or $W$ along the imaginary axis. We have two kinds of estimates: Bounds along a sequence tending to infinity (yielding simple and general bounds for the order of $W_H$ and $W$) and pointwise bounds. Both kinds of bounds have in common that no assumptions are made on the data, except that we ask for limit circle case to hold. We divide our results into two subsections, one for Hamiltonian and one for Jacobi parameters. 

\subsection{Hamiltonian parameters}
\label{S51}

In this subsection we estimate $W_H$ in terms of 
\begin{align}
\label{K68}
b_j^{(s)} \DE \frac{1}{\sqrt{\det \Omega (x_j,x_{j+s})}},
\end{align}
where $s \in \bb N$ can be chosen arbitrarily. The standard choice will be $s=2$, for which
\begin{align}
\label{K67}
b_j^{(2)}=\frac{1}{\sqrt{l_{j+1}l_{j+2}}|\sin (\phi_{j+2}-\phi_{j+1})|}
\end{align}
by \eqref{KX54}. As one can see from \eqref{K68}, the numbers $b_j^{(s)}$ are decreasing in $s$ and hence a larger $s$ leads to a potentially better estimate. 

The following theorem gives a lower bound for $\log |w_{H,22}(ir)|$ for $r$ along an unbounded sequence.

\begin{theorem}
\label{K63}
Let $H$ be a limit circle Hamburger Hamiltonian and choose $s \in \mathbb{N}$, $s \geq 2$. Assume we have two functions $\ms f, \ms g \DF (0,\infty) \to (0,\infty)$ with $\sum_{j=1}^\infty \ms f(b_j^{(s)})=\infty$ and $\sum_{j=1}^\infty \ms g(1/j)<\infty$. Then there exists an increasing and unbounded sequence $(r_m)_{m=1}^\infty$ of positive numbers with
\begin{align*}
\log |w_{H,22}(ir_m)| \gtrsim \frac{1}{\ms g^- (\ms f(r_m))}, \qquad m \in \bb N,
\end{align*}
where $\ms g^-(t) \DE \sup \{r>0 \DF \, \ms g(r) <t\}$.
\end{theorem}

The proof of \Cref{K63} is given further below. First, let us state an important corollary that gives a scale of plain and simple lower bounds for the order of $W_H$. 

\begin{corollary}
\label{K37}
Let $H$ be a limit circle Hamburger Hamiltonian and let $s \in \mathbb{N}$, $s \geq 2$. Then the order of $W_H$ is not less than the convergence exponent of $(b_j^{(s)})_{j =0}^\infty$.
\end{corollary}
\begin{proof}
Let $\alpha$ be the convergence exponent of $(b_j^{(s)})_{j =1}^\infty$. Then, for $\epsilon >0$,
\begin{align*}
\sum_{j=0}^\infty {b_j^{(s)}}^{-(\alpha-\epsilon)} =\infty.
\end{align*}
We can thus apply \Cref{K63} with $\ms f(r)=r^{-(\alpha-\epsilon)}$ and $\ms g(r)=r^{1+\epsilon}$. This results in
\begin{align*}
\log |w_{H,22}(ir_m)| \gtrsim r_m^{\frac{\alpha-\epsilon}{1+\epsilon}}, \qquad m \in \bb N,
\end{align*}
so the order of $W_H$ is at least $(\alpha-\epsilon)/(1+\epsilon)$. Since $\epsilon$ was arbitrary, this finishes the proof.
\end{proof}

\begin{remark}
It is clear from \Cref{K37} that the convergence exponent of $(b_j^{(s)})_{j =0}^\infty$ not larger than $1$. In fact, we have
\begin{align}
\label{K64}
\sum_{j=1}^\infty \frac 1{b_j^{(s)}} \leq s\sqrt{\det \Omega (0,L)} <\infty
\end{align}
which can be seen as follows. For $k \in \{0,\ldots,s-1\}$ we have, due to Minkowski's determinant inequality,
\[
\sqrt{\det \Omega (x_k,L)} \geq \sum_{j=1}^\infty \sqrt{\det \Omega (x_{k+(j-1)s},x_{k+js})}=\sum_{j=1}^\infty \frac{1}{b_{k+(j-1)s}^{(s)}}.
\]
Thus
\[
s\sqrt{\det \Omega (0,L)} \geq \sum_{k=0}^{s-1} \sum_{j=1}^\infty \frac{1}{b_{k+(j-1)s}^{(s)}}=\sum_{j=0}^\infty \frac{1}{b_j^{(s)}}.
\]
\end{remark}

\begin{example}
\label{K93}
Consider a Hamburger Hamiltonian $H$ whose lengths satisfy
\begin{align}
l_j \asymp
\begin{cases}
j^{-\alpha_0}, & j \text{ even}, \\
j^{-\alpha_1}, & j \text{ odd}
\end{cases}
\end{align}
where $\alpha_1>\alpha_0>1$, and whose angles satisfy $\phi_j \DE j \pi/4$. Let us compute the lower bound for the order of $W_H$ from \Cref{K37} for $s=2$ and $s=3$. First we have
\begin{align*}
&\det \Omega (x_j,x_{j+2})=l_{j+1}l_{j+2} \sin^2 (\phi_{j+2}-\phi_{j+1}) \asymp j^{-(\alpha_0+\alpha_1)} \\
&\Rightarrow b_j^{(2)} \asymp j^{\frac{\alpha_0+\alpha_1}{2}}, \qquad \mc E \big((b_j^{(2)})_{j=0}^{\infty} \big)=\frac{2}{\alpha_0+\alpha_1}.
\end{align*}
We turn to $s=3$. Assume first that $j$ is odd. Then
\begin{align*}
&\det \Omega (x_j,x_{j+3})=l_{j+1}l_{j+2} \sin^2 (\phi_{j+2}-\phi_{j+1}) \\
&+l_{j+2}l_{j+3} \sin^2 (\phi_{j+3}-\phi_{j+2})+l_{j+1}l_{j+3} \sin^2 (\phi_{j+3}-\phi_{j+1}) \\
&\asymp (j+1)^{-\alpha_0} (j+2)^{-\alpha_1}+(j+2)^{-\alpha_1} (j+3)^{-\alpha_0}+(j+1)^{-\alpha_0} (j+3)^{-\alpha_0} \\
&\asymp j^{-2\alpha_0}.
\end{align*}
A similar calculation yields $\det \Omega (x_j,x_{j+3}) \asymp j^{-(\alpha_0+\alpha_1)}$ for even $j$. From this we see that 
\begin{align*}
b_j^{(3)} \asymp \begin{cases}
j^{\alpha_0} & j \text{ odd,} \\
j^{\frac{\alpha_0+\alpha_1}{2}} & j \text{ even}
\end{cases}
\end{align*}
and consequently
\begin{align*}
\mc E \big((b_j^{(3)})_{j=0}^{\infty} \big)=\frac{1}{\alpha_0}.
\end{align*}
Considering \Cref{K37}, this shows that $\rho \geq 1/\alpha_0$, where $\rho$ is the order of $W_H$. In fact, also the reverse inequality holds and thus $\rho=1/\alpha_0$; see, e.g., the discussion in the beginning of \Cref{S533} or \Cref{K79}. Note also \Cref{K112}, where we establish that $W_H$ is of positive type with respect to its order.
\end{example}

We shift our focus towards proving \Cref{K63}. The following lemma is needed as preparation.

\begin{lemma}
\label{K80}
Let $(\alpha_j)_{j=1}^\infty$ be a nondecreasing and unbounded sequence of positive numbers. Set
\[
N(r) \DE \# \big\{j \DF \, \alpha_j \leq r   \big\}.
\]
Further, let $\ms f, \ms g \DF (0,\infty) \to (0,\infty)$ be any functions, and assume that $\sum_{j=1}^\infty \ms g(1/j) < \infty$. Then
\begin{align*}
\exists r_0 \, \forall r \geq r_0. \,\, N(r) \leq \frac{1}{\ms g^- (\ms f(r))} \quad \Rightarrow \quad \sum_{j=1}^\infty \ms f(\alpha_j) <\infty,
\end{align*}
where $\ms g^-(t) \DE \sup \{r>0 \DF \, \ms g(r) <t\}$.
\end{lemma}
\begin{proof}
Since $\alpha_1 \leq \alpha_2 \leq \ldots \leq \alpha_{j+1}$ we have $N (\alpha_{j+1}) \geq j+1$. Choosing $j_0$ so large that $\alpha_{j_0+1} \geq r_0$ it follows, for every $j \geq j_0$, that
\begin{align*}
j < N (\alpha_{j+1}) &\leq \frac{1}{\ms g^- (\ms f(\alpha_{j+1}))} \\
&\Rightarrow \ms g^- (\ms f(\alpha_{j+1})) < \frac{1}{j} \\
&\Rightarrow \ms f(\alpha_{j+1}) \leq \ms g \Big( \frac 1j \Big).
\end{align*}
This readily implies
\[
\sum_{j=1}^\infty \ms f(\alpha_j) =\sum_{j=1}^{j_0} \ms f(\alpha_j)+ \sum_{j=j_0}^\infty \ms f(\alpha_{j+1}) \leq  \sum_{j=1}^{j_0} \ms f(\alpha_j)+\sum_{j=j_0}^\infty \ms g(1/j)<\infty.
\]
\end{proof}

\begin{proof}[Proof of \Cref{K63}]
Let $(\alpha_j)_{j=1}^\infty$ be the nondecreasing rearrangement of $(b_j^{(s)})_{j =1}^\infty$, which is an unbounded sequence. By assumption, we have
\[
\sum_{j=1}^{\infty} \ms f(\alpha_j)=\sum_{j=1}^\infty \ms f(b_j^{(s)})=\infty.
\]
Applying \Cref{K80}, we obtain an increasing and unbounded sequence $(r_m)_{m=1}^\infty$ of positive numbers with
\begin{align*}
\#\big\{j \DF \, b_j^{(s)} \leq r_m   \big\}=\#\big\{j \DF \, \alpha_j \leq r_m   \big\}=N(r_m) > \frac{1}{\ms g^- (\ms f(r_m))}, \qquad m \in \bb N.
\end{align*}
Observe that $b_j^{(s)} \leq r_m$ is equivalent to $\sqrt{\det \Omega (x_j,x_{j+s})} \geq 1/r_m$. Since at most $s$ each of the intervals $(x_j,x_{j+s})$ can overlap, for each $m$ we can guarantee the existence of at least
\begin{align*}
k(r_m) \DE \Big\lfloor \frac 1s\#\big\{j \DF \, b_j^{(s)} \leq r_m   \big\} \Big\rfloor
\end{align*}
many disjoint subintervals of $[0,L]$, each of which satisfies $\sqrt{\det \Omega (x_j,x_{j+s})} \geq 1/r_m$. Applying \eqref{K101} once for each $r_m$ gives the lower bound
\begin{align*}
\log |w_{H,22}(ir_m)| \gtrsim k(r_m) \gtrsim \frac{1}{\ms g^- (\ms f(r_m))}.
\end{align*}
\end{proof}

In contrast to \Cref{K63}, which gives a lower bound along a subsequence, the following result gives an estimate for all large enough $r$. It can be seen as a generalisation of \cite[Corollary~2.5]{pruckner.woracek:sinqA}, which only admits $s=2$ and requires a regularity condition. 

\begin{theorem}
\label{K4}
Let $H$ be a limit circle Hamburger Hamiltonian and let $s \in \mathbb{N}$, $s \geq 2$. Suppose that we are given a sequence $(f_j)_{j =1}^\infty$ of nonnegative numbers such that
\[
\sqrt{\det \Omega (x_j,x_{j+s})} \geq f_j, \qquad j \in \mathbb{N}.
\] 
Then 
	\begin{align}
	\label{K78}
	\log |w_{H,22}(ir)| \gtrsim \frac rs \sum_{j=h(r)  }^{\infty} f_j
	\end{align}
	for sufficiently large $r$, where 
\[
h(r) \DE 1+\max \big\{j \in \bb N \DF \, f_j > r^{-1} \big\}.
\]
	The constant implicit in this estimate depends on $H$ but not on $s$.
\end{theorem}

\begin{remark}
\label{K112}
For $f_j=cj^{-\delta}$ (where $\delta>1$) the right-hand side of \eqref{K78} can be explicitly evaluated, and we get
\[
	\log |w_{H,22}(ir)| \gtrsim r^{\frac{1}{\delta}}.
\]
This simple situation occurs, e.g., in \Cref{K93}. Taking $s=4$ we find that
\begin{align*}
\sqrt{\det \Omega (x_j,x_{j+4})} \asymp j^{-\alpha_0}.
\end{align*}
Hence we infer that the growth of $W_H$ is of positive type with respect to its order $1/\alpha_0$ (whereas no type estimate was available in \Cref{K93}). Referring to \cite[Corollary 4.7]{pruckner.reiffenstein.woracek:sinqB-arXiv}, we also see that $W_H$ is of finite type with respect to its order.
%Let us compare this with \Cref{K37}. If $f_j=j^{-\delta}$ is admissible in \Cref{K4}, then $b_j^{(s)} \leq j^{\delta}$ and consequently $\sum_{j=0}^\infty {b_j^{(s)}}^{-(\frac{1}{\delta}}=\infty$. In other words, the convergence exponent of $(b_j^{(s)})_{j=0}^\infty$ is not less than $\frac{1}{\delta}$.
\end{remark}

\begin{proof}[Proof of \Cref{K4}]
Consider an interval $[x_m,x_n]$. Using Minkowski's determinant inequality,
\[
\sqrt{\det \Omega (x_m,x_n)} \geq \sum_{j=1}^{\lfloor \frac{n-m}{s}\rfloor} \sqrt{\det \Omega (x_{m+(j-1)s},x_{m+js})} \geq  \sum_{j=1}^{\lfloor \frac{n-m}{s}\rfloor} f_{m+(j-1)s}.
\]
Applying this estimate also to the intervals $[x_{m+k},x_n]$, we find that
\[
\sqrt{\det \Omega (x_{m+k},x_n)} \geq \sum_{j=1}^{\lfloor \frac{n-m-k}{s}\rfloor} f_{m+k+(j-1)s}.
\]
Summing up the estimates for $k=0,\ldots,s-1$ leads to
\begin{align*}
s&\sqrt{\det \Omega (x_m,x_n)} \geq \sum_{k=0}^{s-1} \sqrt{\det \Omega (x_{m+k},x_n)} \\
&\geq \sum_{k=0}^{s-1} \sum_{j=1}^{\lfloor \frac{n-m-k}{s}\rfloor} f_{m+k+(j-1)s} =\sum_{j=m}^{n-s} f_j.
\end{align*}
Set $m_0(r) \DE h(r)$ and define inductively
\[
m_{n+1}(r) \DE \min \Big\{m>m_n(r) \, : \, \sum_{j=m_n(r)}^m f_j \geq \frac{2s}{r} \Big\}
\]
where $\min \emptyset \DE \infty$. 

Since $f_j \leq 1/r$ for $j \geq h(r)$ we obtain, for each $n$ with $m_{n+1}(r)<\infty$,
\begin{align*}
&\sqrt{\det \Omega (x_{m_n(r)},x_{m_{n+1}(r)})} \geq  \frac 1s \sum_{j=m_n(r)}^{m_{n+1}(r)-s} f_j \\
&\geq \frac 1s \bigg(\sum_{j=m_n(r)}^{m_{n+1}(r)} f_j -\frac sr \bigg) \geq \frac 1r.
\end{align*}
Set $k(r) \DE \min \{n \,:\, m_n(r)=\infty \}$, which is finite due to $\det \Omega (0,L)<\infty$ and Minkowski's determinant inequality. For sufficiently large $r>0$, the ($k(r)+1$ many) points
\[
0, \, x_{m_0(r)}, \, x_{m_1(r)},\ldots x_{m_{k(r)-3}(r)}, \, x_{m_{k(r)-2}(r)}, \, L
\]
meet the condition of \Cref{X67}. Hence, by \eqref{K101},
\[
\log |w_{H,22}(ir)| \gtrsim  k(r)
\]
for large enough $r$. We finish the proof by estimating
\begin{align*}
\frac {2s}r k(r) &\geq  \sum_{n=0}^{k(r)-2} \sum_{j=m_n(r)}^{m_{n+1}(r)-1} f_j + \sum_{j=m_{k(r)-1}(r)}^\infty f_j= \sum_{j=h(r)}^{\infty} f_j.
\end{align*}
\end{proof}

\subsection{Jacobi parameters}
\label{S52}

The numbers $b_j^{(s)}$ defined in \eqref{K68} can be expressed explicitly in terms of the Jacobi parameters $a_n,b_n$. For instance, by comparing \eqref{K62} and \eqref{K67} we see that
\begin{align}
\label{K28}
b_n^{(2)}=\frac{1}{\sqrt{ l_{n+1}l_{n+2}} |\sin (\phi_{n+2} - \phi_{n+1})|}=b_n.
\end{align}
We note that \eqref{K28} together with \eqref{K64} gives an alternative proof of the fact that $\sum_{n=0}^{\infty} b_n^{-1}=\infty$ (Carleman's condition) implies limit point case.

A slightly longer calculation, using \eqref{KX54}, \eqref{K44} and \eqref{K62}, also yields
\begin{align}
\label{K27}
b_n^{(3)} = \frac{b_n b_{n+1}}{\sqrt{a_{n+1}^2+b_n^2+b_{n+1}^2}}.
\end{align}
In this short section we use the identities \eqref{K28}, \eqref{K27} to carry over the lower bounds from \Cref{S51} to the setting of Jacobi matrices. 

For arbitrary $s \geq 2$ we provide more information on how to express $b_j^{(s)}$ in terms of $a_n,b_n$ or in terms of (the values at zero of) orthogonal polynomials in \Cref{K61}. Thus the lower bounds from \Cref{S51} can, in principle, be stated in terms of $a_n,b_n$ even for $s \geq 4$. \newline

We start with the analogue of \Cref{K63}.

\begin{theorem}
\label{K77}
Let $\ms J$ be a Jacobi matrix in limit circle case, and let $W=\smmatrix ACBD$ be its Nevanlinna matrix. Let $\ms f, \ms g \DF (0,\infty) \to (0,\infty)$ with $\sum_{n=1}^\infty \ms g(1/n)<\infty$ and
\begin{align*}
\sum_{n=1}^\infty \ms f(b_n)=\infty \quad \text{ or } \quad \sum_{n=1}^\infty \ms f \bigg(b_n b_{n+1}/\sqrt{a_{n+1}^2+b_n^2+b_{n+1}^2} \bigg)=\infty.
\end{align*}
Then there exists an increasing and unbounded sequence $(r_m)_{m=1}^\infty$ of positive numbers with
\begin{align*}
\log |B(ir_m)| \gtrsim \frac{1}{\ms g^- (\ms f(r_m))}, \qquad m \in \bb N,
\end{align*}
where $\ms g^-(t) \DE \sup \{r>0 \DF \, \ms g(r) <t\}$.
\end{theorem}
\begin{proof}
Let $H$ be the Hamburger Hamiltonian given by \eqref{K11}-\eqref{K62}. We apply \Cref{K63} with $s=2$ and $s=3$, noting \eqref{K28}, \eqref{K27}. Since $w_{H,22}=-B$, this yields the desired lower bound.
\end{proof}

We obtain a remarkably simple lower bound for the order of the Nevanlinna matrix. The slightly weaker estimate $\rho \geq \mc E ((b_n)_{n=0}^\infty )$ has previously been known, but only under (sometimes quite restrictive) regularity conditions. Our result shows, in particular, that no such condition is needed for this estimate to hold.

\begin{corollary} 
\label{K65}
Let $\ms J$ be a Jacobi matrix in limit circle case. Then the order of its Nevanlinna matrix $W$ is not less than the convergence exponent of 
\[
\Big(b_n b_{n+1}/\sqrt{a_{n+1}^2+b_n^2+b_{n+1}^2} \Big)_{n=0}^{\infty}.
\]
In particular, the order of $W$ is not less than the convergence exponent of $(b_n)_{n=0}^{\infty}$.
\end{corollary}
\begin{proof}
Let $H$ be the Hamburger Hamiltonian given by \eqref{K11}-\eqref{K62}.  Since the orders of $W$ and $W_H$ are equal, the theorem follows from \Cref{K37} using \eqref{K28} and \eqref{K27}, and choosing $s=2$ and $s=3$.
\end{proof}

\noindent Using \Cref{K4}, we also get a pointwise lower bound for $\log |B(ir)|$ in terms of the Jacobi parameters.

\begin{theorem}
\label{K59}
Let $\ms J$ be a Jacobi matrix in limit circle case with corresponding Nevanlinna matrix $W=\smmatrix ACBD$. Then
\[
\log |B(ir)| \gtrsim r\sum_{n=h(r)}^{\infty} \frac{1}{b_n}
\]
for sufficiently large $r$, where $h(r) \DE 1+\max \{n \in \bb N \DF \, b_n < r\}$.
\end{theorem}
\begin{proof}
Consider the Hamburger Hamiltonian $H$ given by \eqref{K11}-\eqref{K62}. By \eqref{K28}, \Cref{K4} is applicable with $s=2$ and $f_j =1/b_{j}$. This yields
\[
\log |B(ir)|=\log |w_{H,22}(ir)| \gtrsim r\sum_{n=h(r)}^{\infty} \frac{1}{b_n}.
\]
\end{proof}

\Cref{K59} remains true when $b_n$ is replaced by the term in \eqref{K27}. One proves this in exactly the same way, but taking $s=3$ instead of $s=2$ when applying \Cref{K4}.

\begin{remark}
\label{K61}
For a Hamburger moment problem with orthogonal polynomials $p_n,q_n$, let us express the numbers $b_j^{(s)}$, defined in \eqref{K68} for the associated Hamburger Hamiltonian, in terms of $p_n(0),q_n(0)$. The lengths and angles of the associated Hamburger Hamiltonian are defined recursively by $\phi_1 \DE \pi/2$ and
\begin{align}
&\phi_{j+1}-\phi_j \in (0,\pi), \\
&p_j(0)=\sqrt{l_{j+1}}\sin (\phi_{j+1}), \quad q_j(0)=-\sqrt{l_{j+1}}\cos (\phi_{j+1}).
\end{align}
This is (an index-shifted version of) \cite[(3.21)]{kac:1999}. Clearly this implies
\begin{align*}
\sqrt{l_{j+1}l_{k+1}} \sin (\phi_{j+1}-\phi_{k+1})=q_j(0)p_k(0)-p_j(0)q_k(0).
\end{align*}
Setting $K_{jk} \DE q_j(0)p_k(0)-p_j(0)q_k(0)$, it follows from \eqref{KX54} that
\begin{align}
\label{K58}
\det \Omega (x_m,x_n)= \frac 12 \sum_{j,k=m}^{n-1} K_{jk}^2, \\
\label{K57}
b_n^{(s)}=\det \Omega (x_n,x_{n+s})^{-\frac 12}= \bigg(\frac 12 \sum_{j,k=n}^{n+s-1} K_{jk}^2 \bigg)^{-\frac 12}.
\end{align}
We state two consequences of \eqref{K58}, \eqref{K57}.
\begin{Itemize}
\item In \cite[Remark~5]{braeutigam.mirzoev:2018} the numbers $K_{j+i,j}$ are expressed in terms of $a_n,b_n$, proceeding by induction on $i$. In view of \eqref{K57}, this allows us to write $b_j^{(s)}$ in terms of $a_n,b_n$, in principle for any $s$. However, the expressions become complicated quickly, as can be seen already for $i=1,\ldots,4$:
\begin{align*}
K_{j+1,j} &= \frac{1}{b_j}, \qquad K_{j+2,j}=\frac{a_{n+1}}{b_nb_{n+1}}, \qquad K_{j+3,j} = \frac{b_{j+1}^2+a_{j+1}a_{j+2}}{b_jb_{j+1}b_{j+2}}, \\
 \\
K_{j+4,j} &=\frac{a_{j+1}b_{j+2}^2+a_{j+3}b_{j+1}^2-a_{j+1}a_{j+2}a_{j+3}}{b_jb_{j+1}b_{j+2}b_{j+3}}
.
\end{align*}
\item Relation \eqref{K58} yields the following criterion for indeterminacy: \textit{The Hamburger moment problem with orthogonal polynomials $p_n,q_n$ of the first and second kind is indeterminate if and only if }
\begin{align}
\label{K56}
\sum_{j,k=0}^{\infty} K_{jk}^2=\sum_{j,k=0}^{\infty} \big(q_j(0)p_k(0)-p_j(0)q_k(0) \big)^2<\infty.
\end{align}
This holds because the above sum equals $2 \cdot \det \Omega (0,L)$ (defined for the canonical system associated to the moment problem), and this number is finite if and only if the canonical system is in limit circle case, cf. \cite[Lemma~2.10]{langer.reiffenstein.woracek:kacest-arXiv}.

The indeterminacy criterion \eqref{K56} is arguably simpler than the one given in \cite[Theorem~2]{kostyuchenko.mirzoev:1998}. Surprisingly, the reasoning following \cite[Theorem~2]{kostyuchenko.mirzoev:1998} can be adapted quite easily to obtain another proof for the equivalence of indeterminacy and \eqref{K56}. 
\end{Itemize}
\end{remark}

\section{Upper bounds and situations where the lower bound is attained}
\label{S53}

In this section we estimate $W$ from above, again working with convergence exponents of sequences built from the Hamiltonian or Jacobi parameters. We start by recalling some of what is known for well-behaved data. 

If the lengths as well as the increments of the angles decay like powers (or, more generally, like regularly varying functions) in the sense that
\begin{align}
\label{K95}
l_j \asymp j^{-\alpha }, \qquad |\sin (\phi_{j+1}-\phi_j )| \asymp j^{-\beta},
\end{align}
the lower bound for the order of the Nevanlinna matrix from \Cref{K37} is attained if $\alpha+\beta \geq 2$, i.e., in that case we have $\rho=1/(\alpha +\beta)$ \cite[Theorem 5.3]{pruckner.reiffenstein.woracek:sinqB-arXiv}. The order can be larger than $1/(\alpha +\beta)$ when $\alpha+\beta < 2$, cf. \cite{reiffenstein:kac-hamA-arXiv}.

For data not satisfying \eqref{K95}, a commonly used regularisation method is to estimate the data from above with a power \cite{pruckner.romanov.woracek:jaco}, \cite{pruckner.reiffenstein.woracek:sinqB-arXiv}. A simple example of more roughly behaving data is obtained from \eqref{K95} by changing $l_j$ along the subsequence $j_k=k^2$, demanding $l_{j_k} \asymp j_k^{-\nu /2}= k^{-\nu}$ there instead. For $\nu<2\alpha$ the modified lengths dominate other nearby lengths, so that estimating the lengths from above with a single power would result in a huge loss of information. Yet the peaks occurring in our example are sparse, and their influence on the growth of $W_H$ might be rather small. We will return to this example in \Cref{S534}.

In order to circumvent problems like this we are going to derive upper bounds mostly in terms of convergence exponents, which are more robust to irregularities in the data. In \Cref{S531} we prove an abstract but general upper bound which will be our main tool. Sections \ref{S533} and \ref{S532} are dedicated to one each of two conditions on the angles and the explicit upper estimates that can be obtained under these conditions. We illustrate these results in \Cref{S534}, evaluating the estimates for the example introduced above. Finally, in \Cref{S535} we discuss Berezanskii's classical theorem, which gives sufficient conditions for a Jacobi matrix to be in limit circle case and for the order of $W$ to be equal to the convergence exponent of $(b_n)_{n=0}^{\infty}$. We generalise this theorem by considerably relaxing the conditions on the parameter sequences.

\subsection{A flexible upper bound}
\label{S531}

Our method of obtaining an upper bound for $\log |w_{H,22}(ir)|$ is based on \Cref{K2}. The key in obtaining a good estimate is to efficiently sum up the right-hand sides of \eqref{K24}. The theorem achieves this assuming we have a bound for $\det \Omega$ of a specific form. Our more explicit estimates given in \Cref{S533} and \Cref{S532} are all consequences of \Cref{K26} (and its more abstract version, \Cref{K126}).

\begin{theorem}
\label{K26}
Let $H$ be a limit circle Hamburger Hamiltonian. Suppose that we have two nondecreasing and bounded sequences $(f_n)_{n=0}^\infty,(g_n)_{n=0}^\infty$ of positive numbers, such that
\begin{align}
\label{K23}
\det \Omega (x_m,x_n)^\nu \leq K \cdot \big(f_n-f_m \big)^{\gamma} \big(g_n-g_m\big)^{\delta}, \qquad m<n,
\end{align}
where $\nu,K,\gamma,\delta>0$ and $\gamma+\delta=1$.
Let $f_{\infty} \DE\lim_{n \to \infty} f_n$ and $g_{\infty} \DE \lim_{n \to  \infty} g_n$. 
Then
\begin{align}
\label{K38}
\log |w_{H,22}(ir)| \leq c\frac{K}{\nu} r^{2\nu} (f_{\infty}-f_0)^{\gamma}(g_{\infty}-g_0)^{\delta}
\end{align}
for sufficiently large $r$. The constant $c$ is universal, i.e., it is independent of $H,f_n,g_n,K,\nu,\gamma,\delta$.
\end{theorem}
%\begin{remark}
%The additive term $\BigO (\log (Lr))$ can be neglected at the cost of universality of $c$ and larger $r_0$. To see this it is enough to notice that $w_{H,22}(x_l;.)$ is a polynomial of degree $l$, and hence 
%\begin{align*}
%\log |w_{H,22}(ir)| \geq \log |w_{H,22}(x_l;ir)| \geq (l-1)\log r \geq (l-2)\log (Lr)  
%\end{align*}
%for sufficiently large $r$.
%\end{remark}
\begin{proof}
For fixed $r>0$ we define points $\sigma_k$ and a number $\kappa_H(r) \in \mathbb{N}$ according to the following rules\footnote{This definition is identical with \cite[Definition~5.1]{langer.reiffenstein.woracek:kacest-arXiv}.}.
\begin{itemize}
\item[$\rhd$] Set $\sigma_0 \DE 0$;
\item[$\rhd$] If $\det \Omega (\sigma_{k-1},L) > 1/r^2$, let $\sigma_{k} \in (\sigma_{k-1},L)$ be the unique point satisfying $\det \Omega (\sigma_{k-1},\sigma_{k})=1/r^2$;
\item[$\rhd$] If $\det \Omega (\sigma_{k-1},L) \leq 1/r^2$, set $\sigma_k \DE L$ and $\kappa_H(r) \DE k$, then terminate.
\end{itemize}
Using \Cref{K2} and the estimate $e \log x \leq x^\nu /\nu$ for $x >0$,
\begin{align}
\int_{\sigma_{k-1}}^{\sigma_k} K_H(t;r) \DD t \leq e \log \big(r^2 \det \Omega (\sigma_{k-2},\sigma_k)\big) \\
\leq  \frac {1}{\nu} r^{2\nu} \big( \det \Omega (\sigma_{k-2},\sigma_k) \big)^{\nu}, \quad k=2,\ldots,\kappa_H(r).
\end{align}
We need to pass from $\sigma_k$ to the nearest node $x_n$. Let
\[
n(k) \DE \sup \{n \DF x_n \leq \sigma_k \} \in \bb N \cup \{\infty \}.
\] 
For $k=2,\ldots,\kappa_H(r)-1$, relation \eqref{K23} leads to
\begin{align}
\label{K108}
\begin{split}
&\int_{\sigma_{k-1}}^{\sigma_k} K_H(t;r) \DD t \leq \frac 1\nu r^{2\nu} \big( \det \Omega (x_{n(k-2)},x_{n(k)+1}) \big)^{\nu} \\
&\leq \frac{K}{\nu}  r^{2\nu} \big( f_{n(k)+1}-f_{n(k-2)} \big)^{\gamma}\big( g_{n(k)+1}-g_{n(k-2)} \big)^{\delta},
\end{split}
\end{align}
and this is still true for $k=\kappa_H(r)$, interpreting $x_{\infty +1}, \, f_{\infty +1}, \, g_{\infty +1}$ as $L,f_{\infty},g_{\infty}$, respectively.
 We take the sum and apply H{\"o}lder's inequality to obtain
\begin{align}
\label{K109}
\begin{split}
&\int_{\sigma_1}^{L} K_H(t;r) \DD t \\
&\leq \frac{K}{\nu} r^{2\nu} \sum_{k=2}^{\kappa_H(r)} \big( f_{n(k)+1}-f_{n(k-2)} \big)^{\gamma}\big( g_{n(k)+1}-g_{n(k-2)} \big)^{\delta} \\
&\leq \frac{K}{\nu} r^{2\nu} \Big(\sum_{k=2}^{\kappa_H(r)}  f_{n(k)+1}-f_{n(k-2)} \Big)^{\gamma} \Big(\sum_{k=2}^{\kappa_H(r)}  g_{n(k)+1}-g_{n(k-2)} \Big)^{\delta}.
\end{split}
\end{align}
Note that $n(k)+1 \leq n(k+1)$ for $k=2,\ldots,\kappa_H(r)-1$. Both of these sums are thus estimated from above by a telescoping sum:
\begin{align*}
\sum_{k=2}^{\kappa_H(r)} f_{n(k)+1} &-f_{n(k-2)} \leq f_{\infty}-f_{n(\kappa_H(r)-2)}\\
&+\sum_{i=1}^3 \sum_{\substack{k=2 \\ k \equiv i \text{ mod } 3}}^{\kappa_H(r)-1} \big( f_{n(k+1)}-f_{n(k-2)} \big) \leq 3 \big( f_{\infty}-f_{0}\big)  
\end{align*}
and analogously for $(g_n)$ in place of $(f_n)$. In total, we get
\begin{align}
\label{K97}
&\int_{\sigma_1}^{L} K_H(t;r) \DD t \leq \frac{3K}{\nu} r^{2\nu} \big(f_{\infty}-f_0\big)^{\gamma} \big(g_{\infty}-g_0\big)^{\delta}.
\end{align}
We need to estimate the contribution of $[0,\sigma_1]$. Assuming $r>1/ \sqrt{\det \Omega (0,L)}$, \cite[Proposition~5.8]{langer.reiffenstein.woracek:kacest-arXiv} gives
\begin{align}
&\int_0^{\sigma_1} K_H(t;r) \DD t \leq 2\log (Lr)+1-\log 4 \leq 2\log (Lr).
\end{align}
Invoking \eqref{X65}, we obtain
\begin{align}
\log |w_{H,22}(ir)| \leq c'\Big(\frac{3K}{\nu} r^{2\nu} (f_{\infty}-f_0)^{\gamma}(g_{\infty}-g_0)^{\delta}+2\log (Lr)\Big)
\end{align}
for $r>1/ \sqrt{\det \Omega (0,L)}$, where $c'$ is a universal constant. For sufficiently large $r$ we can incorporate the term $2\log (Lr)$ into the term left of the plus sign, at the cost of increasing $c'$ by a small margin. This finishes the proof.
\end{proof}

In \cite[Corollary~5.3]{pruckner.woracek:sinqA} Hamiltonians of the form $H(t)=\xi_{\phi (t)} \xi_{\phi (t)}^{\top}$ with continuously rotating angle $\phi (t)$ are considered. It is shown that for H\"older continuous $\phi$ with exponent $\alpha \in (0,1]$, the monodromy matrix is of order at most $1/(1+\alpha)$ with finite type. The following result shows that a similar bound holds for sequences of angles satisfying a discrete H\"older-type condition.

\begin{corollary}
Let $H$ be a limit circle Hamburger Hamiltonian. If there are positive numbers $\alpha,d$ such that
\[
|\phi_{m+1}-\phi_n| \leq d|x_m-x_n|^\alpha, \qquad m < n,
\]
then
\begin{align}
\label{K48}
\log |w_{H,22}(ir)| \leq r^{\frac{1}{1+\alpha}} \cdot 2cL(1+\alpha)d^{\frac{1}{1+\alpha}}
\end{align}
for large enough $r$, where $c$ is the universal constant from \Cref{K26}.
\end{corollary}
\begin{proof}
Let $m<n$. Then, using \eqref{KX54},
\begin{align}
\label{K36}
\begin{split}
\det \Omega (x_m,x_n) &= \frac 12 \sum_{j,k=m+1}^n l_jl_k \sin^2(\phi_j-\phi_k) \\
&\leq \frac{d^2}2 (x_n-x_m)^{2\alpha} \sum_{j,k=m+1}^n l_jl_k \leq d^2 (x_n-x_m)^{2\alpha +2}.
\end{split}
\end{align}
Hence \eqref{K23} is satisfied with $f_n \DE g_n \DE x_n$, exponents $\nu=1/(2(1+\alpha))$ and $\gamma=\delta=1/2$, and constant $K=d^{1/(1+\alpha)}$. The assertion follows from \Cref{K26} applied with these parameters.
\end{proof}

\Cref{K26} takes into account the behaviour of $H$ on the whole interval $[0,L]$. It is well-suited to prove Theorems \ref{K89} and \ref{K79} below, which are based on $\ell^p$-conditions on the parameters. With little additional effort we can prove the following extension which is needed in the proofs of \Cref{K66} and \Cref{K49}. It estimates the integral in \eqref{X65} on a given subinterval of $(0,L)$, and allows more than two sequences to help estimating $\det \Omega$.

\begin{proposition}
\label{K126}
Let $H$ be a limit circle Hamburger Hamiltonian and let $p \in \bb N$. Suppose that we have
\begin{itemize}
\item[$\rhd$] numbers $\gamma_i>0$, $i=1,\ldots,p$, with $\sum_{i=1}^p \gamma_i=1$;
\item[$\rhd$] $M \in \bb N \cup \{0\}$, $N \in \bb N \cup \{\infty\}$, with $M<N$;
\item[$\rhd$] nondecreasing and bounded sequences $(f_n^{(i)})_{n=M}^N$, $i=1,\ldots, p,$ of positive numbers (If $N=\infty$, we interpret $f_{\infty}^{(i)}$ as the limit of $f_n^{(i)}$ as $n \to \infty$);
\item[$\rhd$] $\nu,K>0$
\end{itemize}
such that
\begin{align}
\label{K123}
\det \Omega (x_m,x_n)^\nu \leq K \cdot \prod_{i=1}^p\big(f_n^{(i)}-f_m^{(i)} \big)^{\gamma_i}, \qquad M \leq m<n \leq N.
\end{align}
Then
\begin{align}
\label{K138}
\int_{x_M}^{x_N} K_H(t;r) \DD t \leq \frac{3K}{\nu} r^{2\nu} \prod_{i=1}^p \big(f_N^{(i)}-f_M^{(i)}\big)^{\gamma_i}+2(1+e)\log (Lr)
\end{align}
for $r>1/ \sqrt{\det \Omega (0,L)}$ (where $x_{\infty} \DE L$).
\end{proposition}
\begin{proof}
We adapt the proof of \Cref{K26} in the following way. For fixed $r>0$ we define $\sigma_k$ only in the interval $[x_M,x_N]$:
\begin{itemize}
\item[$\rhd$] Set $\sigma_0 \DE x_M$;
\item[$\rhd$] If $\det \Omega (\sigma_{k-1},x_N) > 1/r^2$, let $\sigma_{k} \in (\sigma_{k-1},x_N)$ be the unique point satisfying $\det \Omega (\sigma_{k-1},\sigma_{k})=1/r^2$;
\item[$\rhd$] If $\det \Omega (\sigma_{k-1},x_N) \leq 1/r^2$, set $\sigma_k \DE x_N$, $d \DE k$, and terminate.
\end{itemize}
As in the proof of \Cref{K26} we define
\[
n(k) \DE \sup \{n \DF x_n \leq \sigma_k \} \in \bb N \cup \{\infty \}.
\] 
Let $k \in \{2,\ldots,d-1 \}$. Starting with the same estimate as in the upper row of \eqref{K108}, we now use \eqref{K123} to estimate $\det \Omega$ from above:
\begin{align}
\label{K125}
&\int_{\sigma_{k-1}}^{\sigma_k} K_H(t;r) \DD t \leq \frac{K}{\nu}  r^{2\nu} \prod_{i=1}^p \big( f_{n(k)+1}^{(i)}-f_{n(k-2)}^{(i)} \big)^{\gamma_i}.
\end{align}
In a similar manner, for $k=d$ we obtain
\begin{align*}
&\int_{\sigma_{d-1}}^{\sigma_d} K_H(t;r) \DD t \leq \frac{K}{\nu}  r^{2\nu} \prod_{i=1}^p \big( f_N^{(i)}-f_{n(d-2)}^{(i)} \big)^{\gamma_i}
\end{align*}
which, upon interpreting $f_{n(d)+1}^{(i)}$ as $f_N^{(i)}$, is of the same form as \eqref{K125}.
Summing up as in \eqref{K109}, \eqref{K97}, recalling $\sigma_d=x_N$ and using H{\"o}lder's inequality for multiple sequences, yields
\begin{align}
\label{K98}
&\int_{\sigma_1}^{x_N} K_H(t;r) \DD t \leq \frac{K}{\nu} r^{2\nu} \sum_{k=2}^{d} \prod_{i=1}^p \big( f_{n(k)+1}^{(i)}-f_{n(k-2)}^{(i)} \big)^{\gamma_i} \\
\nonumber
&\leq \frac{K}{\nu} r^{2\nu} \prod_{i=1}^p \Big(\sum_{k=2}^{d}  f_{n(k)+1}^{(i)}-f_{n(k-2)}^{(i)} \Big)^{\gamma_i} \leq \frac{3K}{\nu} r^{2\nu} \prod_{i=1}^p \big(  f_{N}^{(i)}-f_{M}^{(i)} \big)^{\gamma_i}.
\end{align}
In order to estimate the contribution of $[x_M,\sigma_1]$ we use \cite[Proposition~5.8]{langer.reiffenstein.woracek:kacest-arXiv} for $[x_M,\hat t(r)]$ and \Cref{K2} for $[\hat t(r),\sigma_1]$. Assuming $r>1/ \sqrt{\det \Omega (0,L)}$, we have
\begin{align}
\label{K99}
\begin{split}
&\int_{x_M}^{\hat t(r)} K_H(t;r) \DD t \leq 2 \log (Lr), \\
&\int_{\hat t(r)}^{\sigma_1} K_H(t;r) \DD t \leq e \log \big(r^2 \det \Omega (0,\sigma_1)\big) \leq 2e \log (Lr)
\end{split}
\end{align}
which, combined with \eqref{K98}, proves \eqref{K138}.
\end{proof}

\subsection{Exploiting summability of angle increments}
\label{S533}

In \cite[Theorem~4.7]{berg.szwarc:2014} it is shown that, if the orthogonal polynomials $p_n,q_n$ of a Hamburger moment sequence are such that
\begin{align}
\label{K106}
(p_n(0)^2)_{n=0}^{\infty}, (q_n(0)^2)_{n=0}^{\infty} \in \ell^p
\end{align}
for some $p>0$, then the order of $W$ is not larger than $p$ (and $W$ is of finite type with respect to $p$). Since the lengths of the associated Hamburger Hamiltonian are given by
\begin{align}
\label{K107}
l_{n+1}=p_n(0)^2+q_n(0)^2, \qquad n \in \bb N \cup \{0\},
\end{align}
cf. \cite[(3.22)]{kac:1999}, the condition \eqref{K106} is equivalent to $(l_j)_{j=1}^{\infty} \in \ell^p$. This fact was noticed in \cite[Proposition 2.3]{pruckner.romanov.woracek:jaco}.

We give new upper bounds that also take into account the angles. In this subsection we consider angles whose step sizes (modulo $\pi$) are summable, and prove an upper bound for the order of the Nevanlinna matrix that is typically correct, e.g., if the data satisfies \eqref{K95}. For angles whose step sizes are not summable we refer the reader to \Cref{S532}.

\begin{theorem}
\label{K89}
Let $H$ be a limit circle Hamburger Hamiltonian with lengths $(l_j)_{j=1}^\infty$ and angles $(\phi_j)_{j=1}^\infty$. Suppose we have $\alpha, \beta \geq 1$ such that
\begin{align*}
(l_j)_{j=1}^\infty \in \ell^{\frac{1}{\alpha}}, \qquad (|\sin (\phi_{j+1}-\phi_j)|))_{j=1}^\infty \in \ell^{\frac{1}{\beta}}.
\end{align*}
Then
\begin{align}
\label{K87}
\begin{split}
\log &|w_{H,22}(ir)| \\
&\leq r^{\frac{1}{\alpha+\beta}} \cdot 2c (\alpha+\beta) \Big(\sum_{j=1}^{\infty} l_j^{\frac 1{\alpha}} \Big)^{\frac{\alpha}{\alpha+\beta}}\Big(1+\sum_{j=1}^{\infty} |\sin(\phi_{j+1}-\phi_j)|^{\frac 1{\beta}} \Big)^{\frac{\beta}{\alpha+\beta}}
\end{split}
\end{align}
for large enough $r$, where $c$ is the same constant as in \Cref{K26}.
\end{theorem}

Before we turn to the proof of this theorem, let us state an immediate corollary which gives an upper bound just for the order.

\begin{corollary}
\label{K104}
Let $H$ be a limit circle Hamburger Hamiltonian and assume that $\sum_{j=1}^{\infty} |\sin (\phi_{j+1}-\phi_j)| <\infty$. Set 
\begin{align*}
\alpha_0 &\DE \frac 1{\mc E \big(\big(\frac{1}{l_j} \big)_{j=1}^{\infty} \big)} = \sup \Big\{\alpha>0 \, \big| \, \sum_{j=1}^{\infty} l_j^{\frac 1\alpha}<\infty \Big\}, \\
\beta_0 &\DE \frac 1{\mc E \big(\big(\frac{1}{|\sin (\phi_{j+1}-\phi_j)|} \big)_{j=1}^{\infty} \big)} = \sup \Big\{\beta>0 \, \big| \, \sum_{j=1}^{\infty} |\sin(\phi_{j+1}-\phi_j)|^{\frac 1\beta}<\infty \Big\}.
\end{align*}
Then
\begin{align*}
\rho \leq \frac{1}{\alpha_0+\beta_0}.
\end{align*}
\end{corollary}

\begin{proof}[{Proof of \Cref{K89}}]
We repeatedly use the inequality $|\sin (a+b)| \leq |\sin a|+|\sin b|$ for $a,b \in \mathbb{R}$, which is easy to check. With help of \eqref{KX54} we estimate
\begin{align}
\label{K100}
\begin{split}
\det \Omega (x_m,x_n)&= \sum_{\substack{j,k=m+1 \\ j<k}}^n l_jl_k \sin^2(\phi_j-\phi_k) \\
&\leq \sum_{\substack{j,k=m+1 \\ j<k}}^n l_jl_k \Big(\sum_{s=j}^{k-1} |\sin (\phi_{s+1}-\phi_s)| \Big)^2 \\
&\leq   \Big(\sum_{j=m+1}^n l_j \Big)^2 \Big(\sum_{j=m+1}^n |\sin (\phi_{j+1}-\phi_j)| \Big)^2.
\end{split}
\end{align}
We have
\begin{align}
\label{K96}
\Big(\sum_{j=m+1}^n l_j \Big)^{\frac 1\alpha} \leq \sum_{j=m+1}^n l_j^{\frac 1\alpha}
\end{align}
due to the triangle inequality in $\ell^{\alpha}$. An analogous estimate holds for $|\sin (\phi_{j+1}-\phi_j)|$ in place of $l_j$ and $\beta$ in place of $\alpha$. Hence
\begin{align*}
\det \Omega (x_m,x_n)\leq  \Big(\sum_{j=m+1}^n l_j^{\frac{1}{\alpha}} \Big)^{2\alpha} \Big(\sum_{j=m+1}^n |\sin (\phi_{j+1}-\phi_j)|^{\frac{1}{\beta}} \Big)^{2\beta},
\end{align*}
i.e., \eqref{K23} is satisfied with $f_n=\sum_{j=1}^n l_j^{\frac 1{\alpha}}$, $g_n=\sum_{j=1}^n |\sin(\phi_{j+1}-\phi_j)|^{\frac 1{\beta}}$, exponents $\nu=1/(2(\alpha+\beta))$, $\gamma=\alpha/(\alpha+\beta)$, $\delta=\beta/(\alpha+\beta)$, and constant $K=1$. The statement follows from \Cref{K26}. 
\end{proof}

If the suprema defining $\alpha_0,\beta_0$ in \Cref{K104} are maxima, we expect \Cref{K89} to give a good upper estimate for $\log |w_{H,22}(ir)|$ itself. Otherwise the theorem does not give a finer estimate than just the order. In order to say a bit more about $\log |w_{H,22}(ir)|$ in this case we can use the following result, provided we know the decay of tails of the sums of lengths and angle differences.

\begin{proposition}
\label{K66}
Let $H$ be a limit circle Hamburger Hamiltonian and suppose that $\sum_{j=1}^\infty |\sin (\phi_{j+1}-\phi_j )| <\infty $. Let $G$ be the decreasing function 
\begin{align*}
G \FD{\mathbb{N}}{(0,\infty)}{1}{N}{ \frac 1N \big(\sum_{j=N+1}^\infty l_j \big)^{\frac 12} \big( \sum_{j=N+1}^\infty |\sin (\phi_{j+1}-\phi_j)| \big)^{\frac 12} }
\end{align*}
Then
	\[
	\log |w_{H,22}(ir)| \lesssim G^- \Big(\frac{\log r}{\sqrt r} \Big) \cdot \log r
	\]
	for large enough $r$, where $G^-(x)=\min \{N \in \bb N \DF \, G(N) < x \}$.
\end{proposition}
\begin{proof}
It follows from \eqref{K100} that \eqref{K123} is satisfied with sequences $f_n^{(1)} \DE x_n$, $f_n^{(2)} \DE \sum_{j=1}^n |\sin (\phi_{j+1}-\phi_j)|$, exponents $\nu=1/4$, $\gamma_1=\gamma_2=1/2$, and constant $K=1$. \Cref{K126} gives, for arbitrary $N \in \bb N \cup \{0\}$,
\begin{align}
&\int_{x_N}^L  K_H(t;r) \DD t \leq 12 r^{\frac 12} \cdot NG(N)+2(1+e)\log (Lr).
\end{align}
We estimate the integral over $[x_0,x_N]$ by separately applying \Cref{K126} on every interval $[x_{j-1},x_j]$, $j=1,\ldots,p$. Since $\det \Omega (x_{j-1},x_j)=0$, the condition \eqref{K123} is satisfied for any choice of $p,f_n^{(i)},\gamma_i,\nu,K$. Letting $K$ tend to zero we obtain
\begin{align*}
\int_{x_{j-1}}^{x_j} K_H(t;r) \DD t \leq 2(1+e)\log (Lr).
\end{align*}
Summing up we obtain
\begin{align*}
\int_{0}^{x_N} K_H(t;r) \DD t \leq 2(1+e)\cdot N\log (Lr).
\end{align*}
The upper bounds on $[0,x_N]$ and $[x_N,L]$ become balanced when setting $N \DE G^-  \big(\log r /\sqrt r \big)$:
\begin{align}
\label{K102}
\begin{split}
\int_{0}^L  K_H(t;r) \DD t &\leq 12 G^- \Big(\frac{\log r}{\sqrt r}\Big)\log r+ \BigO(\log r) \\&+2(1+e) G^- \Big(\frac{\log r}{\sqrt r}\Big)\log r \big(1+\BigO \big((\log r)^{-1}\big) \big).
\end{split}
\end{align}
We finish the proof by incorporating the terms within $\BigO(.)$ into the dominating term.
\end{proof}

\subsection{Exploiting convergence of angles}
\label{S532}

When the angles are not of bounded variation, we expect the growth of the Nevanlinna matrix to depend on whether or not the angles are convergent and, if they do converge, on their speed of convergence (this expectation is based on what happens for regularly varying data, cf. \cite{reiffenstein:kac-hamA-arXiv}). Indeed this intuition is reflected in the following result, which is applicable to any limit circle Hamburger Hamiltonian. The decay of the lengths and the convergence of the angles are measured by $\ell^p$-conditions which, as mentioned earlier, are less sensible to perturbations than previously used conditions.

\begin{theorem}
\label{K79}
Let $H$ be a limit circle Hamburger Hamiltonian with lengths $(l_j)_{j=1}^\infty$ and angles $(\phi_j)_{j=1}^\infty$. Suppose we have $\psi \in \bb R$ and $\alpha, \omega \geq 1$ such that
\begin{align*}
(l_j)_{j=1}^\infty \in \ell^{\frac{1}{\alpha}}, \qquad (l_j \sin^2 (\phi_j-\psi))_{j=1}^\infty \in \ell^{\frac{1}{\omega}}.
\end{align*}
Then
\begin{align}
\label{K47}
\begin{split}
\log &|w_{H,22}(ir)| \\
&\leq r^{\frac{2}{\alpha+\omega}} \cdot c (\alpha+\omega) \Big(\sum_{j=1}^{\infty} l_j^{\frac 1{\alpha}} \Big)^{\frac{\alpha}{\alpha+\omega}}\Big(\sum_{j=1}^{\infty} \big(l_j \sin^2 (\phi_j-\psi)\big)^{\frac 1{\omega}} \Big)^{\frac{\omega}{\alpha+\omega}}
\end{split}
\end{align}
for sufficiently large $r$, where $c$ is the universal constant from \Cref{K26}.
\end{theorem}

We prove the theorem after stating the resulting upper bound for the order:

\begin{corollary}
\label{K94}
Let $H$ be a limit circle Hamburger Hamiltonian and fix $\psi \in \bb R$. Set 
\begin{align*}
\alpha_0 &\DE \frac 1{\mc E \big(\big(\frac{1}{l_j} \big)_{j=1}^{\infty} \big)} = \sup \Big\{\alpha>0 \, \big| \, \sum_{j=1}^{\infty} l_j^{\frac 1\alpha}<\infty \Big\}, \\
\omega_0 &\DE \frac 1{\mc E \big(\big(\frac{1}{l_j\sin^2 (\phi_j-\psi)} \big)_{j=1}^{\infty} \big)} = \sup \Big\{\omega>0 \, \big| \, \sum_{j=1}^{\infty} \big(l_j \sin^2(\phi_j-\psi)\big)^{\frac 1\omega}<\infty \Big\}.
\end{align*}
Then
\begin{align*}
\rho \leq \frac{2}{\alpha_0+\omega_0}.
\end{align*}
\end{corollary}

%\begin{example}[{Continuation of \Cref{K93}}]
%\label{K96}
%We would like to treat the parameters from \Cref{K93} with \Cref{K94}. In order to get a good upper bound, we assume that the angles converge as fast as possible (under the condition \eqref{K103} on the angle increments), i.e., there exists $\psi \in \bb R$ such that
%\begin{align*}
%|\sin (\phi_j-\psi)| \leq j^{-\beta}.
%\end{align*}
%We then have
%\begin{align*}
%\alpha_0 = \min \{\nu,\alpha \}, \qquad \omega_0=\min \{\nu+4\beta,\alpha+2\beta \}.
%\end{align*}
%Hence a sufficient condition for $\rho=\frac{1}{\alpha+\beta}$ is that $\alpha \leq \nu$. We do not know the precise value of $\rho$ when $\nu<\alpha$. Nonetheless, we have the bounds
%\begin{align*}
%\frac{1}{\alpha+\beta} \leq \rho \leq \max\Big\{\frac{1}{\nu+2\beta},\frac{1}{\frac{\nu+\alpha}{2}+\beta} \Big\},
%\end{align*}
%whereas the upper bound obtained from \cite[Corollary~2.5]{pruckner.reiffenstein.woracek:sinqB-arXiv} is $2/(\nu+2\beta)$ and thus worse.
%\end{example}

\begin{proof}[{Proof of \Cref{K79}}]
For integers $n>m$ and $\psi \in \bb R$ we estimate $\det \Omega$ by
\begin{align}
\label{K92}
\begin{split}
&\det \Omega (x_m,x_n) = \det \Big(\sum_{j=m+1}^n l_j \xi_{\phi_j} \xi_{\phi_j}^\top \Big) \\
&= \det \bigg[\begin{pmatrix}
\cos \psi & \sin \psi \\
-\sin \psi & \cos \psi
\end{pmatrix}
\sum_{j=m+1}^n l_j\xi_{\phi_j} \xi_{\phi_j}^\top 
 \begin{pmatrix}
\cos \psi & -\sin \psi \\
\sin \psi & \cos \psi
\end{pmatrix} \bigg] \\
&= \det  \Big(\sum_{j=m+1}^n l_j\xi_{\phi_j-\psi} \,\xi_{\phi_j-\psi}^\top \Big) \\
&\leq \Big(\sum_{j=m+1}^n l_j \cos^2 (\phi_j-\psi) \Big) \Big( \sum_{j=m+1}^n l_j \sin^2 (\phi_j-\psi) \Big) \\
&\leq \Big(\sum_{j=m+1}^n l_j \Big) \Big( \sum_{j=m+1}^n l_j \sin^2 (\phi_j-\psi) \Big).
\end{split}
\end{align}
We further estimate this term using the triangle inequalities in $\ell^{\alpha}$ and $\ell^{\omega}$, arriving at
\begin{align*}
\det \Omega (x_m,x_n) \leq \Big(\sum_{j=m+1}^n l_j^{\frac 1\alpha} \Big)^{\alpha} \Big( \sum_{j=m+1}^n \big(l_j \sin^2 (\phi_j-\psi)\big)^{\frac 1\omega} \Big)^{\omega}.
\end{align*}
The statement follows from an application of \Cref{K26} with parameters $f_n=\sum_{j=1}^n l_j^{\frac 1{\alpha}}$, $g_n=\sum_{j=1}^n \big(l_j \sin^2(\phi_j-\psi)\big)^{\frac 1{\omega}}$, exponents $\nu=1/(\alpha+\omega)$, $\gamma=\alpha/(\alpha+\omega)$, $\delta=\omega/(\alpha+\omega)$, and constant $K=1$.
\end{proof}

The upper bound in \Cref{K79} becomes better the faster the angles converge. Yet it does not factor in (the decay of) the differences of consecutive angles, which is an important piece of information. In the theorem below this information is also taken account of, in addition to the potential convergence of the angles.

\begin{theorem}
\label{K49}
Let $H$ be a limit circle Hamburger Hamiltonian. Assume we have $\alpha \geq 1$ and $\beta \in (0,1)$ such that
\begin{align*}
(l_j)_{j=1}^\infty \in \ell^{\frac{1}{\alpha}}, \qquad (|\sin (\phi_{j+1}-\phi_j)|))_{j=1}^\infty \in \ell^{\frac{1}{\beta}}.
\end{align*} 
Fix $\psi \in \mathbb{R}$ and consider the increasing function 
\begin{align*}
F \FD{\mathbb{N}}{(0,\infty)}{1}{N}{ N^{\frac{1-\beta}{\alpha}} \big[\big(\sum_{j=N+1}^\infty l_j \big) \big( \sum_{j=N+1}^\infty l_j \sin^2 (\phi_j-\psi) \big)\big]^{-\frac {\alpha+1}{2\alpha}}. }
\end{align*}
Then
	\[
	\log |w_{H,22}(ir)| \lesssim \big(rF^-(r)^{1-\beta} \big)^{\frac{1}{\alpha+1}}
	\]
	for large enough $r$, where $F^-(r)=\min \{N \in \bb N \DF \, F(N) < r \}$.
\end{theorem}

We postpone the proof to the end of the subsection and continue with an immediate corollary that gives an estimate for the order. We point out that the decay of $|\sin (\phi_{j+1}-\phi_j)|$ influences the upper bound for the order, even if the angles do not converge.

\begin{corollary}
\label{K91}
Let $H$ be a limit circle Hamburger Hamiltonian and set
\begin{align}
\label{K72}
\alpha_0 &\DE \frac 1{\mc E \big(\big(\frac{1}{l_j} \big)_{j=1}^{\infty} \big)} = \sup \Big\{\alpha>0 \, \big| \, \sum_{j=1}^{\infty} l_j^{\frac 1\alpha}<\infty \Big\}, \\
\label{K75}
\beta_0 &\DE \frac 1{\mc E \big(\big(\frac{1}{|\sin (\phi_{j+1}-\phi_j)|} \big)_{j=1}^{\infty} \big)} = \sup \Big\{\beta>0 \, \big| \, \sum_{j=1}^{\infty} |\sin(\phi_{j+1}-\phi_j)|^{\frac 1\beta}<\infty \Big\}.
\end{align}
Assume $\alpha_0>1$ and $\beta_0 \in (0,1]$, and that $\sum_{j=N+1}^\infty l_j \lesssim N^{1-(\alpha_0-\epsilon )}$ for every small\footnote{as is the case when $l_j \lesssim j^{-\alpha_0}$, even for $\epsilon=0$.} $\epsilon>0$. Then
\begin{align}
\label{K121}
\rho \leq \frac{\alpha_0-\beta_0}{\alpha_0^2-\beta_0}.
\end{align}
\end{corollary}
\begin{proof}
Fixing a small $\epsilon>0$ we apply \Cref{K49} with $\wt{\alpha} \DE \alpha_0-\epsilon$ and $\wt{\beta} \DE \beta_0-\epsilon$ and estimating $\sin^2(\phi_j-\psi)$ with $1$. Then
\begin{align*}
F(N) &\geq N^{\frac{\wt{\alpha}^2-\wt{\beta}}{\wt{\alpha}}}
\end{align*}
and consequently
\begin{align*}
\log |w_{H,22}(ir)| \lesssim r^{\frac{\wt{\alpha}-\wt{\beta}}{\wt{\alpha}^2-\wt{\beta}}}.
\end{align*}
Letting $\epsilon$ go to zero we obtain the asserted upper bound for the order.
\end{proof}

\begin{remark}
If $\beta_0=1$ then the upper bound in \eqref{K121} may be attained, for if $l_j \asymp j^{-\alpha_0}$ and $|\sin (\phi_{j+1}-\phi_j)| \asymp j^{-1}$ the lower bound from \Cref{K37} is $1/(\alpha_0+1)$. We do not know whether \eqref{K121} is sharp for other values of $\beta$, but we have reason to believe that it is not. This is because, in the specific situation of \eqref{K95}, a better estimate is available in the literature; for such lengths and angles we have $\alpha_0=\alpha$ and $\beta_0=\beta$, and \cite[Example 2.23]{pruckner.romanov.woracek:jaco} states
\begin{align}
\label{K74}
\rho \leq \begin{cases}
\frac{1}{\alpha_0+\beta_0} & \text{if } \alpha_0+\beta_0 \geq 2 \\
\frac{1-\beta_0}{\alpha_0-\beta_0} & \text{if } \alpha_0+\beta_0 < 2.
\end{cases}
\end{align}
We have seen in \Cref{K89} that summability of the angle increments is a sufficient condition for \eqref{K74}, if $\alpha_0,\beta_0$ are defined by \eqref{K72} and \eqref{K75} (note that we have $\beta_0 \geq 1$ if the angle increments are summable).
\end{remark}

\begin{GenericTheorem}{Open Question}
Let $H$ be a limit circle Hamburger Hamiltonian and define $\alpha_0,\beta_0$ by \eqref{K72}, \eqref{K75}. Does \eqref{K74} always hold, even if $\sum_{j=1}^{\infty} |\sin (\phi_{j+1}-\phi_j)|=\infty$?
\end{GenericTheorem}

\begin{proof}[{Proof of \Cref{K49}}]
For each $r$ we cut $[0,L]$ into two pieces $[0,x_{T(r)}]$ and $[x_{T(r)},L]$ where $T(r)$ is to be determined later. 
\item[\textit{Contribution of $[0,x_{T(r)}]$}:] 
By \eqref{K100},
\begin{align*}
\begin{split}
\det \Omega (x_m,x_n)\leq  \Big(\sum_{j=m+1}^n l_j \Big)^2 \Big(\sum_{j=m+1}^n |\sin (\phi_{j+1}-\phi_j)| \Big)^2.
\end{split}
\end{align*}
Using \eqref{K96} and applying H{\"o}lder's inequality to the sum over $1 \cdot |\sin (\phi_{j+1}-\phi_j)|$ we obtain
\begin{align*}
\det \Omega (x_m,x_n) &\leq   \Big(\sum_{j=m+1}^n l_j^{\frac 1\alpha} \Big)^{2\alpha}  \Big(\sum_{j=m+1}^n |\sin (\phi_{j+1}-\phi_j)|^{\frac 1\beta} \Big)^{2\beta} (n-m)^{2-2\beta}.
\end{align*}
Hence \eqref{K123} is satisfied for $\nu=1/(2(\alpha+1))$, $K=1$ and
\begin{align*}
f_n^{(1)} &\DE \sum_{j=1}^n l_j^{\frac 1\alpha},& f_n^{(2)} &\DE \sum_{j=1}^n |\sin (\phi_{j+1}-\phi_j)|^{\frac 1\beta}, &f_n^{(3)} &\DE n, \\
\gamma_1 &\DE \frac{\alpha}{\alpha+1}, &\gamma_2 &\DE \frac{\beta}{\alpha+1}, &\gamma_3 &\DE \frac{1-\beta}{\alpha+1}.
\end{align*}
\Cref{K126} gives
\begin{align*}
&\int_0^{x_{T(r)}} K_H(t;r) \DD t \\
&\leq 6(\alpha+1)r^{\frac{1}{\alpha+1}} \Big(\sum_{j=1}^{T(r)} l_j^{\frac{1}{\alpha}} \Big)^{\frac{\alpha}{\alpha+1}} \Big(\sum_{j=1}^{T(r)} |\sin (\phi_{j+1}-\phi_j)|^{\frac 1\beta} \Big)^{\frac{\beta}{\alpha+1}} T(r)^{\frac{1-\beta}{\alpha+1}} \\
&+2(1+e)\log (Lr).
\end{align*}
In view of the assumptions we made on summability, the estimate simplifies to
\begin{align}
\label{K111}
&\int_0^{x_{T(r)}} K_H(t;r) \DD t \lesssim \big( rT(r)^{1-\beta}\big)^{\frac{1}{\alpha+1}}+\BigO\big(\log (Lr) \big). 
\end{align}
\item[\textit{Contribution of $[x_{T(r)},L]$}:] By \eqref{K92},
\begin{align*}
&\det \Omega (x_m,x_n) \leq \Big(\sum_{j=m+1}^n l_j \Big) \Big( \sum_{j=m+1}^n l_j \sin^2 (\phi_j-\psi) \Big).
\end{align*}
This shows that \eqref{K123} is satisfied with parameters $f_n^{(1)} \DE x_n$, $f_n^{(2)} \DE  \sum_{j=1}^n l_j \sin^2 (\phi_j-\psi )$, exponents $\nu=\gamma_1=\gamma_2=1/2$, and constant $K=1$. Hence \Cref{K126} yields
\begin{align}
\nonumber
&\int_{x_{T(r)}}^L  K_H(t;r) \DD t \\
\nonumber
&\leq 6r \cdot \Big(\sum_{j=T(r)+1}^{\infty} l_j \Big)^{\frac 12} \Big( \sum_{j=T(r)+1}^{\infty} l_j \sin^2 (\phi_j-\psi) \Big)^{\frac 12}+2(1+e)\log (Lr) \\
\label{K114}
&\asymp r T(r)^{\frac{1-\beta}{\alpha+1}} F(T(r))^{-\frac{\alpha}{\alpha+1}}+\BigO\big(\log (Lr) \big).
\end{align}
We equate this to \eqref{K111} (neglecting the $\BigO(.)$-terms) and arrive at $F(T(r))=r$. We finish the proof by plugging $T(r) \DE F^-(r)$ into \eqref{K111} and noting that the resulting term is an upper bound for both parts of the integral, as $F(F^-(r)) < r$ and $\log (Lr) \lesssim r^{1/(\alpha+1)}$.
\end{proof}

\subsection{Evaluating the upper bounds for a concrete example}
\label{S534}

Let us illustrate the applicability of the results from the previous two sections by returning to the example sketched in the introduction to \Cref{S53}. Namely, consider a Hamburger Hamiltonian $H$ whose lengths satisfy
\begin{align}
\label{K117}
l_j \asymp
\begin{cases}
j^{-\frac{\nu}2}, & j=k^2 \text{ for some } k \in \bb N, \\
j^{-\alpha}, & \text{else}
\end{cases}
\end{align}
with $\nu,\alpha>1$. Further assume that the angles $\phi_j$ of $H$ satisfy
\begin{align}
\label{K116}
|\sin (\phi_{j+1}-\phi_j)| \asymp j^{-\beta}
\end{align}
for $\beta \geq 0$. The universal lower bound from \Cref{K37} with $s=2$ takes the form $\rho \geq 1/(\alpha+\beta)$ (as a straightforward calculation shows). Below we discuss the upper bounds that we obtain from Theorems \ref{K89}, \ref{K79} and \ref{K49}, some of which coincide with the universal lower bound. We distinguish between four cases. \\

\noindent \textbf{Case 1: $\beta>1$.} This is the most well-behaved case. With notation as in \Cref{K104} we have
\begin{align*}
\alpha_0 = \min \{\alpha,\nu \}, \qquad \beta_0=\beta.
\end{align*}
Hence
\[
\rho \leq \frac{1}{\min \{\alpha,\nu \}+\beta}
\]
by \Cref{K104}. This implies $\rho=1/(\alpha+\beta)$ when $\nu \geq \alpha$. We do not know the exact value of $\rho$ when $\nu <\alpha$. \\

\noindent \textbf{Case 2: $\beta \leq 1$, angles converge as fast as possible.}
Under the condition \eqref{K116}, the fastest possible decay of $|\sin (\phi_j-\psi)|$ is at the rate of $j^{-\beta}$, in the sense that $|\sin (\phi_j-\psi)| \lesssim j^{-\gamma}$ implies $\gamma \leq \beta$. 
%This follows simply from the triangle inequality:
%\[
%j^{-\beta} \asymp |\sin (\phi_{j+1}-\phi_j)| \leq |\sin (\phi_{j+1}-\psi)|+|\sin (\phi_j-\psi)|.
%\]
So, let us assume that \eqref{K117} holds and that $|\sin (\phi_j-\psi)| \lesssim j^{-\beta}$. Then, with notation from \Cref{K94},
\begin{align*}
\alpha_0 = \min \{\nu,\alpha \}, \qquad \omega_0=\min \{\nu+4\beta,\alpha+2\beta \}.
\end{align*}
Hence 
\begin{align*}
\rho \leq
\begin{cases}
\frac{1}{\alpha+\beta} & \nu \geq \alpha \\
\frac{1}{\frac{\nu+\alpha}{2}+\beta} & \alpha-2\beta \leq \nu < \alpha \\
\frac{1}{\nu+2\beta} & \nu < \alpha-2\beta.
\end{cases}
\end{align*}
We see that, for $\nu \geq \alpha$, the upper bound coincides with the universal lower bound and thus $\rho=1/(\alpha+\beta)$. We do not know the precise value of $\rho$ when $\nu<\alpha$. Nonetheless, our upper bound is better than the one obtained from \cite[Corollary~2.5]{pruckner.reiffenstein.woracek:sinqB-arXiv}, which (for $\nu <\alpha$) is $2/(2\beta+\nu+1)$.\\

\noindent \textbf{Case 3: $\beta \leq 1$, no convergence of angles.} In this case we 
cannot determine $\rho$ (except in special cases), but we can give the following upper bound which becomes better as $\beta$ increases:
\begin{align}
\label{K120}
\rho \leq \begin{cases}
\frac{\alpha-\beta}{\alpha^2-\beta} & \nu \geq 2\alpha-1 \\[.5ex]
\frac{\nu+1-2\beta}{(\nu-1)(\alpha+1)+2-2\beta} & \alpha \leq \nu < 2\alpha-1 \\[.5ex]
\frac{\nu+1-2\beta}{\nu^2+1-2\beta} & \nu <\alpha.
\end{cases}
\end{align}
This upper bound for the order is different from the universal lower bound $1/(\alpha+\beta)$, except if $\nu \geq \alpha$ and $\beta=1$. Note that for $\nu \geq 2\alpha$ we do not lose much information when estimating $l_j$ with $j^{-\alpha}$. Thus we expect  \cite[Corollary 2.5]{pruckner.reiffenstein.woracek:sinqB-arXiv} to give a good estimate. Indeed we obtain a better estimate than \eqref{K120} in that case:
\begin{align*}
\nu \geq 2\alpha \,  \Rightarrow \,\rho \, \begin{cases}
\leq \frac{1-\beta}{\alpha-\beta} & \alpha+\beta<2 \\[.5ex]
=\frac{1}{\alpha+\beta} & \alpha+\beta \geq 2.
\end{cases}
\end{align*}
\begin{proof}[{Proof of \eqref{K120}}]
We are going to apply \Cref{K49}. For a small $\epsilon>0$ and $\wt{\alpha} \DE \min\{\alpha,\nu\}-\epsilon$, $\wt{\beta} \DE \beta - \epsilon<1$ we have that
\begin{align*}
(l_j)_{j=1}^\infty \in \ell^{\frac{1}{\wt{\alpha}}}, \qquad (|\sin (\phi_{j+1}-\phi_j)|))_{j=1}^\infty \in \ell^{\frac{1}{\wt{\beta }}}.
\end{align*}
We also need to choose $\psi \in \bb R$, but this choice does not matter since we are going to estimate $\sin^2 (\phi_j-\psi) $ with $1$. Since
\begin{align}
\label{K124}
&\sum_{\substack{j=N+1 \\ j \neq k^2}}^{\infty} l_j \asymp \sum_{j=N+1}^{\infty} j^{-\alpha} \asymp N^{1-\alpha}, \qquad \sum_{\substack{j=N+1 \\ j=k^2}}^{\infty} l_j \asymp \sum_{k=\lceil \sqrt{N+1} \rceil}^{\infty} k^{-\nu} \asymp N^{\frac{1-\nu}{2}}
\end{align}
we have
\begin{align*}
F(N) &\geq N^{\frac{1-\wt{\beta}}{\wt{\alpha}}}\Big(\sum_{j=N+1}^{\infty} l_j \Big)^{-\frac{\wt{\alpha}+1}{\wt{\alpha}}} \asymp N^{\frac{1-\wt{\beta}}{\wt{\alpha}}} \big(N^{\max \{1-\alpha,\frac{1-\nu}{2}\}} \big)^{-\frac{\wt{\alpha}+1}{\wt{\alpha}}} \\
&\Rightarrow F^-(r) \lesssim r^{\frac{\wt{\alpha}}{1-\wt{\beta}+\min \{\alpha-1,\frac{\nu-1}{2}\} (\wt{\alpha}+1)}}.
\end{align*}
The upper bound for $\log |w_{H,22}(ir)|$ given by \Cref{K49} is thus not larger than a constant times $r^{\tau}$, where
\begin{align}
\label{K115}
\begin{split}
\tau &= \frac{1-\wt{\beta}+\min \{\alpha-1,\frac{\nu-1}{2}\} (\wt{\alpha}+1)+\wt{\alpha}(1-\wt{\beta})}{(1-\wt{\beta}+\min \{\alpha-1,\frac{\nu-1}{2}\} (\wt{\alpha}+1))(\wt{\alpha}+1)} \\
&=\frac{1-\wt{\beta}+\min \{\alpha-1,\frac{\nu-1}{2}\}}{1-\wt{\beta}+\min \{\alpha-1,\frac{\nu-1}{2}\} (\wt{\alpha}+1)}.
\end{split}
\end{align}
In particular we see that $\rho \leq \tau$, and this holds true for every small $\epsilon>0$. Hence we still get an upper bound for $\rho$ if on the right-hand side of \eqref{K115} we replace $\wt{\alpha}$ by $\min\{\alpha,\nu \}$ and $\wt{\beta}$ by $\beta$, resulting in \eqref{K120}.
\end{proof}

\noindent \textbf{Case 4 (mixed): $\beta \leq 1$, angles converge like a given power.}
In this case we assume that the angles converge at an intermediate rate, i.e., $|\sin (\phi_j-\psi)| \lesssim j^{-\gamma}$ with $\gamma \leq \beta$, in addition to \eqref{K117}, \eqref{K116}. For $\gamma=\beta$ we are in the situation of\footnote{The reason we kept cases $2$ and $3$ separate is to emphasise the extreme cases: fast convergence versus no convergence of the angles.} case $2$, while $\gamma=0$ corresponds to case $3$. The upper bound for the order obtained from \Cref{K49} is
\begin{align}
\label{K118}
\rho \leq \begin{cases}
\frac{\alpha-\beta+\gamma}{\alpha^2-\beta+(\alpha+1)\gamma} & \nu \geq 2\alpha-1 \\[.5ex]
\frac{\nu+1-2\beta+2\gamma}{(\nu-1)(\alpha+1)+2-2\beta+2(\alpha+1)\gamma} & \alpha \leq \nu < 2\alpha-1 \\[.5ex]
\frac{\nu+1-2\beta+2\gamma}{\nu^2+1-2\beta+2(\nu+1)\gamma} & \nu <\alpha.
\end{cases}
\end{align}
We note that for $\gamma=0$ this is the same upper bound as in case $3$ (because we use \Cref{K49} once more). For $\gamma=\beta$ we are in case $2$ and one can check that the upper bound obtained there is either the same as or better than the one in \eqref{K118}. 
\begin{proof}[{Proof of \eqref{K118}}]
First, we have
\begin{align*}
&\sum_{\substack{j=N+1 \\ j \neq k^2}}^{\infty} l_j\sin^2 (\phi_j-\psi) \lesssim \sum_{j=N+1}^{\infty} j^{-\alpha-2\gamma} \asymp N^{1-\alpha-2\gamma}, \\
&\sum_{\substack{j=N+1 \\ j=k^2}}^{\infty} l_j\sin^2 (\phi_j-\psi) \lesssim \sum_{k=\lceil \sqrt{N+1} \rceil}^{\infty} k^{-\nu-4\gamma} \asymp N^{\frac{1-\nu-4\gamma}{2}}.
\end{align*}
As in case $2$ we set $\wt{\alpha} \DE \min\{\alpha,\nu\}-\epsilon$, $\wt{\beta} \DE \beta - \epsilon<1$ for some small $\epsilon>0$. Thus, recalling \eqref{K124}, we can estimate the functions $F$ and $F^-$ from \Cref{K49} as follows:
\begin{align*}
F(N) &\gtrsim N^{\frac{1-\wt{\beta}}{\wt{\alpha}}}\Big(N^{\max \{1-\alpha,\frac{1-\nu}{2}\}-\gamma} \Big)^{-\frac{\wt{\alpha}+1}{\wt{\alpha}}} \\
&\Rightarrow F^-(r) \lesssim r^{\frac{\wt{\alpha}}{1-\wt{\beta}+(\wt{\alpha}+1) [\min \{\alpha-1,\frac{\nu-1}{2}\}+\gamma ] }}.
\end{align*}
\Cref{K49} yields
\begin{align}
\log |w_{H,22}(ir)| &\lesssim \big(rF^-(r)^{1-\wt{\beta}} \big)^{\frac{1}{\wt{\alpha}+1}}=r^{\tau}, \\
\tau &=\frac{1-\wt{\beta}+\min \{\alpha-1,\frac{\nu-1}{2}\}+\gamma}{1-\wt{\beta}+(\wt{\alpha}+1) [\min \{\alpha-1,\frac{\nu-1}{2}\}+\gamma] }.
\end{align}
Letting $\epsilon$ go to zero we arrive at \eqref{K118}.
\end{proof}

%But what happens when we replace $\phi_j$ with
%\begin{align*}
%\widehat{\phi}_j \DE j^{-\beta}, 
%\end{align*}
%i.e., the angles converge monotonically to their limit instead of alternatingly? Let $\widehat{\rho}$ be the order of the Hamburger Hamiltonian with lengths $l_j$ and angles $\widehat{\phi}_j$. Analogously, assuming $\alpha \leq \nu$ the upper bound we get from \Cref{K94} is
%\begin{align*}
%\widehat{\rho} \leq \frac{1}{\alpha+\beta}.
%\end{align*}
%However, we have
%\begin{align*}
%|\widehat{\phi}_{j+1}-\widehat{\phi}_j|=j^{-\beta}-(j+1)^{-\beta} \sim \beta j^{-\beta-1}.
%\end{align*}
%\Cref{K37} gives the lower bound
%\begin{align}
%\label{K105}
%\widehat{\rho} \geq \frac{1}{\alpha+\beta+1}
%\end{align}
%which is actually attained, due to \Cref{K104}. Here it is the step size of the angles that matters, not their rate of convergence. Compare this to yet another sequence of angles,
%\begin{align*}
%\widetilde{\phi}_j \DE (-1)^j j^{-\beta-1},
%\end{align*}
%which converges faster to $0$ but with essentially the same step size, namely
%\begin{align*}
%|\widetilde{\phi}_{j+1}-\widetilde{\phi}_j| \asymp j^{-\beta-1}.
%\end{align*}
%The order of the Hamiltonian with lengths $l_j$ and angles $\widetilde{\phi}_j$ again equals $1/(\alpha+\beta+1)$.

\subsection{Generalising Berezanskii's theorem}
\label{S535}

We return to the setting of Jacobi matrices and consider again the lower bound for the order of the Nevanlinna matrix given by the convergence exponent of $(b_n)_{n=0}^\infty$, cf. \Cref{K65}. Although the order can be strictly larger than that, cf. \cite[Theorem~1.3]{reiffenstein:kac-hamA-arXiv}, equality holds in a surprising number of previously studied situations. Examples are Berezanskii's theorem \cite[Theorems~4.2 and 4.11]{berg.szwarc:2014} and Jacobi parameters with power asymptotics in the subcritical case \cite[Theorem~1]{pruckner:blubb}. These situations have in common that all solutions $(u_n)_{n=0}^\infty$ of the recurrence
\begin{align}
\label{K82}
b_n u_{n+1} + a_n u_n+b_{n-1}u_{n-1} =0, \qquad n \geq 1
\end{align}
satisfy $u_n=\BigO (b_{n-1}^{-\frac 12})$. This property, together with Carleman's condition being violated, is actually all we need for the order to be equal to the convergence exponent of $(b_n)_{n=0}^\infty$. In contrast to all of our previous results, prevalence of the limit circle case is not a prerequisite, but a consequence of the other assumptions.

\begin{proposition}
\label{K84}
Let $\ms J$ be a Jacobi matrix for which Carleman's condition is violated, i.e., $\sum_{n=0}^\infty b_n^{-1}<\infty$. If all solutions $(u_n)_{n=0}^\infty$ of \eqref{K82} satisfy
\begin{align}
\label{K86}
u_n=\BigO (b_{n+k}^{-\frac 12}), \qquad n \to \infty
\end{align}
for a fixed integer $k$, then $\ms J$ is in limit circle case and the order $\rho$ of its Nevanlinna matrix $W$ is equal to the convergence exponent of $(b_n)_{n=0}^\infty$. If, additionally,
\[
\sum_{n=0}^\infty \frac{1}{b_n^{\rho}}<\infty
\]
then $W$ is of finite type with respect to its order.
\end{proposition}
\begin{proof}
By \eqref{K107} and \eqref{K86} we have $l_{n+1}=\BigO(b_n^{-1})$. Hence the prevalence of limit circle case follows from the assumption $\sum_{n=0}^\infty b_n^{-1}<\infty$.
In addition, $(b_n^{-1})_{n=0}^{\infty} \in \ell^p$ implies $(l_j)_{j=1}^{\infty} \in \ell^p$ and thus $\rho \leq p$, either by \cite[Theorem~4.7]{berg.szwarc:2014} or by \Cref{K79}. Clearly we also get finite type with respect to $p$, and with respect to $\rho$ itself provided that $(b_n^{-1})_{n=0}^{\infty} \in \ell^{\rho}$. The lower bound for $\rho$ is given by \Cref{K65}.
\end{proof}

Under the assumptions of \cite[Theorem~4.2]{berg.szwarc:2014}, decay of solutions as in \eqref{K86} is a consequence of direct manipulations of the recurrence relation \eqref{K82}. Still, the assumptions can be relaxed considerably: the diagonal parameters $a_n$ need not be $\Smallo (.)$ of the off-diagonal parameters $b_n$, and oscillation of a periodic type can be admitted, cf. \cite{swiderski.trojan:2017} (easier to read) and \cite{swiderski.trojan:2020} (more general). Another reference is \cite[Theorem~3.9]{yafaev:2020}, which is significantly less general. However, the assumptions made there can be more directly compared to those of \cite{berg.szwarc:2014}, and they lead us to the following Berezanskii-type theorem. 

\begin{theorem}
\label{K85}
Let $\ms J$ be a Jacobi matrix with diagonal parameters $a_n \in \bb R$ and off-diagonal parameters $b_n>0$. Let 
\[
\beta_n \DE -\frac{a_n}{2\sqrt{b_{n-1}b_n}}.
\]
Assume
\begin{Enumerate}
\item $\sum_{n=1}^\infty b_n^{-1}<\infty$ (Carleman's condition violated);
\item $\sum_{n=1}^\infty |\beta_{n+1}-\beta_n|<\infty$ and $\lim_{n \to \infty} \beta_n \in (-1,1)$ (relative smallness and regularity of diagonal);
\item Regularity of the off-diagonal parameters:
\[
\sum_{n=1}^\infty \Big|\frac{b_n}{\sqrt{b_{n-1}b_{n+1}}}-1\Big|<\infty.
\]
\end{Enumerate}
Then $\ms J$ is in limit circle case and the order $\rho$ of its Nevanlinna matrix $W$ is equal to the convergence exponent of $(b_n)_{n=0}^\infty$. If, additionally,
\[
\sum_{n=0}^\infty \frac{1}{b_n^{\rho}}<\infty
\]
then $W$ is of finite type with respect to its order.
\end{theorem}
\begin{proof}
It follows from \cite[Theorem~3.9]{yafaev:2020} that \eqref{K82} has a solution $(f_n)_{n=0}^\infty$ satisfying
\[
f_n=\frac{1}{\sqrt{b_n}}e^{-i\theta_n} (1+\Smallo (1))
\]
where $\theta_n=\sum_{m=1}^{n-1} \arccos \beta_m$ and the sum is taken over all $m$ with $|\beta_m| \leq 1$. We see that $\overline{f_n}$ is a linearly independent solution of the same recursion, and hence we are in the situation of \Cref{K84}.
\end{proof}

\begin{remark}
\Cref{K85} is a significant generalisation of Berezanskii's theorem. In order to see this we compare the assumptions of \Cref{K85} to those of Berezanskii's theorem, as stated in \cite[Theorem~4.11]{berg.szwarc:2014}. Apart from the violated Carleman condition the classical assumptions are $\sum_{n=1}^\infty |\beta_n|<\infty$ (which is much more restrictive than the second condition in \Cref{K85}) and log-concavity or log-convexity of $(b_n)_{n=0}^\infty$:
\[
			\Big(\exists n_0\in\bb N\DQ\forall n\geq n_0\DP b_n^2\geq b_{n-1}b_{n+1}\Big)
			\text{ or }
			\Big(\exists n_0\in\bb N\DQ\forall n\geq n_0\DP b_n^2\leq b_{n-1}b_{n+1}\Big).
		\]
The perhaps more interesting case is the log-concave one, since in the log-convex case the parameters $b_n$ are always exponentially growing (leading to $\rho=0$).
 Let us show that log-concavity together with summability of the sequence $(1/b_n)_{n=0}^{\infty}$ implies the third condition in \Cref{K85}. The proof is similar to that of \cite[Lemma~4.1]{berg.szwarc:2014}. \\
We notice first that log-concavity implies that $b_{n+1}/b_n$ is decreasing and thus converges to a limit $\lambda$. In addition, $\lambda \geq 1$ since otherwise Carleman's condition would be satisfied. Hence, for $n_0 \in \bb N$ sufficiently large,
\begin{align*}
&\sum_{n=n_0}^\infty \bigg|\frac{b_n}{\sqrt{b_{n-1}b_{n+1}}}-1\bigg| =\sum_{n=n_0}^\infty \sqrt{\frac{b_n}{b_{n-1}}}\Bigg|\sqrt{\frac{b_n}{b_{n+1}}}-\sqrt{\frac{b_{n-1}}{b_n}}\Bigg| \\
&\leq \sqrt{\lambda+1} \sum_{n=n_0}^\infty \Bigg|\sqrt{\frac{b_n}{b_{n+1}}}-\sqrt{\frac{b_{n-1}}{b_n}}\Bigg|=\sqrt{\lambda+1} \Bigg(\frac{1}{\sqrt{\lambda}}-\sqrt{\frac{b_{n_0-1}}{b_{n_0}}} \Bigg)<\infty.
\end{align*}
\end{remark}

%---------
%   FINISH
%---------

\printbibliography

{\footnotesize

\begin{flushleft}
	J.~Reiffenstein \\
	Department of Mathematics\\
	Stockholms universitet\\
	106 91 Stockholm\\
	SWEDEN\\
	email: jakob.reiffenstein@math.su.se \\[5mm]
\end{flushleft}

}

\end{document}